\documentclass[reqno,12pt]{amsart}
\usepackage{amsfonts,amsmath,amsthm,amssymb,amscd}
\usepackage{latexsym}

\catcode`\@=11

\@addtoreset{equation}{section} \@addtoreset{equation}{subsection}
\def\theequation{\ifnum\value{subsection}>0\relax
\thesubsection.\arabic{equation}\relax
\else\ifnum\value{section}>0\relax
\thesection.\arabic{equation}\relax \else\arabic{equation}\fi\fi}

\newtheorem{thm}[equation]{Theorem}
\newtheorem{lemma}[equation]{Lemma}

\newtheorem{cor}[equation]{Corollary}

\newfont {\gothic}{eufm10 at 10pt}
\keywords{Classical invariant theory, coinvariants, degenerate principal series, Howe correspondence}

\begin{document}

\title [Degenerate principal series and Howe correspondence] {Degenerate principal series of metaplectic groups and Howe correspondence}

\author [] {Soo Teck Lee}
\address{Department of Mathematics\\
National University of Singapore\\
10 Lower Kent Ridge Road\\
Singapore 119076}
\email{matleest@nus.edu.sg}

\author [] {Chen-Bo Zhu}
\address{Department of Mathematics\\
National University of Singapore\\
10 Lower Kent Ridge Road\\
Singapore 119076}
\email{matzhucb@nus.edu.sg}

\date{}
\thanks{Both authors are supported by the MOE2010-T2-2-113.}
\maketitle

\begin{abstract} The main purpose of this article is to supplement the authors' results on degenerate principal series representations of real symplectic groups with the analogous results for metaplectic groups. The basic theme, as in the previous case, is that their structures are anticipated by certain natural subrepresentations constructed from Howe correspondence. This supplement is necessary as these representations play a key role in understanding the basic structure of Howe correspondence (and its complications in the archimedean case), and their global counterparts play an equally essential part in the proof of Siegel-Weil formula and its generalizations (work of Kudla-Rallis). The full results in the metaplectic case also shed light on the seeming peculiarities, when the results in the symplectic case are viewed in their isolation.

\end{abstract}

\section{Introduction: classical invariant theory and its transcendental analog}
\label{Intro}

Let the complex orthogonal group ${\mathrm{O}_m}$ act on the space of complex matrices $M_{m,n}$ by matrix multiplication on the left:
\[{\mathrm{O}_m}\curvearrowright M_{m,n}.\]

The First Fundamental Theorem (FFT) of classical invariant theory asserts that the ring of ${\mathrm{O}_m}$-invariant polynomials are generated by the fundamental invariants of degree $2$:
\[P(M_{m,n})^{{\mathrm{O}_m}}=<r_{ij}|\ 1\leq i,j\leq n>,\]
where $r_{ij}(X)=\sum _{k=1}^mx_{ki}x_{kj}$, for $X=(x_{ij})\in M_{m,n}$.

Let $S^2({\mathbb C}^n)$ denote the space of complex $n\times n$ symmetric matrices and define
\[Q:\ M_{m,n} \rightarrow S^2({\mathbb C}^n),\]
\[X\mapsto X^tX. \]
Then FFT asserts that the pull-back map
\[Q^*:\ P(S^2({\mathbb C}^n)) \rightarrow P(M_{m,n})^{{\mathrm{O}_m}}\]
is surjective. Alternatively, we have an affine embedding:
\[\psi:\ M_{m,n}//{\mathrm{O}_m} \hookrightarrow S^2({\mathbb C}^n),\]
where $M_{m,n}//{\mathrm{O}_m}$ is the affine quotient of $M_{m,n}$ by ${\mathrm{O}_m}$. The Second Fundamental Theorem (SFT) of classical invariant theory is then a statement on the image of $\psi$ as an affine subvariety of $S^2({\mathbb C}^n)$. Note that both $M_{m,n}$ and $S^2({\mathbb C}^n)$ carry natural actions of $\mathrm{GL}(n,{\mathbb C})$, and all maps ($Q, Q^*, \psi$) are $\mathrm{GL}(n,{\mathbb C})$-equivariant.

\vspace{0.2in}

We now state results of Kudla-Rallis \cite{KR1} and Lee-Zhu \cite{LZ2}, which may be viewed as transcendental analogs of FFT and SFT alluded to above.

Let \[{\mathrm{O}(p,q)}\curvearrowright M_{p+q,n}({\mathbb R}), \]
again by matrix multiplication on the left. Denote by $G={\widetilde{\mathrm{Sp}}(2n,{\mathbb R})}$ the real metaplectic group of rank $n$ and $P$ its Siegel parabolic subgroup. The rest of notation will be explained in Section \ref{Embed}.

%{\bf Theorem A}:

\begin{thm} {\rm (Kudla-Rallis)} There is a natural topological embedding with closed image:
\[\psi _{p,q}:{\mathcal S}(M_{p+q,n}({\mathbb R}))_{{\mathrm{O}(p,q)}}\hookrightarrow {I^\alpha(\sigma)} (:=\mbox{\rm Ind}_{P}^{G}(\chi ^{\alpha}_{\sigma})),\]
as $G$-representations.
\end{thm}

%{\bf Theorem B}:

\begin{thm} {\rm (Lee-Zhu)} describes the image of $\psi _{p,q}$.
\end{thm}

In \cite{LZ2}, the authors describe the image $\Omega ^{p,q}$ of $\psi _{p,q}$ when $p+q$ is even, in which case the representations concerned factor through the linear group ${\mathrm{Sp}(2n,{\mathbb R})}$. In fact only the subcase of $p-q\equiv 0$ (mod $4$) was treated in detail and the other subcase of $p-q\equiv 2$ (mod $4$) was left to the reader.

The aim of the current article is to complete the description of $\Omega ^{p,q}$ (without any restriction on $p$ and $q$). This is necessary and useful as the representations involved (the coinvariants) play a key role in understanding the basic structure of Howe correspondence (through the doubling method \cite{H1,Ra,Ku}) and their global counterparts feature prominently in the proof of Siegel-Weil formula and its generalizations \cite{KR3}. For a recent application of these results to first occurrence conjecture of Kudla-Rallis, see \cite{SZ}. As an additional benefit, this full description enables us to organize the statements in a more coherent way so that their structures emerge clearer to the reader. The basic idea is very simple: since both $\Omega ^{p,q}$ and ${I^\alpha(\sigma)} $ are $K$-multiplicity free ($K$ being a maximal compact subgroup of $G$), and since the structure of ${I^\alpha(\sigma)}$ is known (Section \ref{Subquo}; \cite{L2} for the linear group), one can identify the image $\Omega ^{p,q}$ as a $G$-representation by knowing its $K$-types. The fascinating point is that from this description, one concludes that the reducibilities of $I^\alpha(\sigma)$ are completely accounted for by the possible embeddings of $\Omega ^{p,q}$'s (in a precise way). We summarize this assertion in Theorem \ref{summary}, which should be viewed as analogous to SFT (proverbially stated) as ``all relations among the fundamental invariants are generated by the obvious ones".

Here are some words on the organization of this article. In Section \ref{Howe}, we review the basics of Howe duality correspondence. In Section \ref{Embed}, we introduce the space of coinvariants, as a special case of Howe's construction of maximal quotients, and its embedding into the degenerate principal series $I^{\alpha}(\sigma)$ of the metaplectic group. In Section \ref{TC}, we describe the transition coefficients of $I^{\alpha}(\sigma)$ and give its immediate consequences for irreducibility and complementary series. In Section \ref{Subquo}, we give the detailed structure of $I^{\alpha}(\sigma)$ at points of reducibility, which are again arrived by analyzing the transition coefficients. In Section \ref{Image}, we complete our description of the image of $\psi _{p,q}$ in $I^{\alpha}(\sigma)$. As proofs of (any new) results follow a similar line as those of \cite{L2} or \cite{LZ2}, we shall omit them.

Finally we mention the following works which are closely related to the theme of this article: \cite{KR1,LZ1,LZ2,LZ3,Ya} (for real groups) and
\cite{KR2,KS,Ya} (for $p$-adic groups).

\section{Howe duality correspondence}
\label{Howe}

In this section, we review briefly Howe's theory of dual pair correspondence, for the case at hand.

Let $H={\mathrm{O}(p,q)}$, which acts on $M_{p+q,n}({\mathbb R})$ by matrix multiplication on the left. We have the dualised action of $H$ on ${\mathcal Y} = L^2(M_{p+q,n}({\mathbb R}))$, the space of square integrable functions on $M_{p+q,n}({\mathbb R})$. As is well-known, there is a (unitary) representation $\omega$ of a metaplectic group $\widetilde{\rm{Sp}}(2N,{\mathbb R})$ ($N=(p+q)n$) on ${\mathcal Y}$, called the Schrodinger model of an oscillator representation. Here and after, for any real symplectic group ${\rm Sp}(2N,{\mathbb R})$, $\widetilde{\rm{Sp}}(2N,{\mathbb R})$ denotes its metaplectic two fold cover. The oscillator representation depends on a choice of a non-trivial unitary character of ${\mathbb R}$, which we fix once for all.

The main feature of this set-up is the following: there is a reductive dual pair $$(H,G)=({\mathrm{O}(p,q)}, {\mathrm{Sp}(2n,{\mathbb R})})\subseteq {\rm Sp}(2N,{\mathbb R}),$$
namely a pair of reductive subgroups in ${\rm Sp}(2N,{\mathbb R})$ which are mutual centralizers.

Let $\omega ^{\infty}$ be the smooth representation of $\omega$ realized in ${\mathcal Y}^{\infty}$. For the case at hand, ${\mathcal Y}^{\infty}={\mathcal S}(M_{p+q,n}({\mathbb R}))$, the Schwartz space of rapidly decreasing functions on $M_{p+q,n}({\mathbb R})$.

For any subgroup $E$ of ${\rm Sp}(2N,{\mathbb R})$, denote by $\widetilde{E}$ the preimage of $E$ under the covering map $\widetilde{{\rm Sp}}(2N,{\mathbb R})\rightarrow {\rm Sp}(2N,{\mathbb R})$. If $E$ is also reductive, denote by ${\mathcal R} (\widetilde{E}, \omega)$ the set of infinitesimal equivalent classes of irreducible admissible representations of $\widetilde{E}$ which are realizable as quotients by $\omega ^{\infty}(\widetilde{E})$-invariant closed subspaces of ${\mathcal Y}^{\infty}$.

\vspace{0.2in}

\noindent {\bf Howe Duality Theorem} {\rm (\cite{H2})}: ${\mathcal R}(\widetilde{H}\cdot\widetilde{G},\omega)$ is the graph of a bijection between ${\mathcal R}(\widetilde{H},\omega)$ and
${\mathcal R}(\widetilde{G},\omega)$.

\vspace{0.2in}

This is the Howe quotient correspondence. Let us make it more
concrete. As usual, we may twist the oscillator representation (by a
character of $\widetilde{H}$) so that it will factor through the
standard linear action of $H$ on ${\mathcal Y}$. We will assume that
this has been done (but retain the same notation) and will be
concerned with $H$, rather than $\widetilde{H}$. Let $\rho $ be an
irreducible admissible representation of $H$, which is realizable as
a quotient by a $H$-invariant closed subspace of ${\mathcal
Y}^{\infty}$. Define $\Omega (\rho )$ to be the maximal quotient of
${\mathcal Y}^{\infty}$ on which $H$ acts by a representation of
class $\rho $. We have (isomorphism class to mean infinitesimal
equivalent class):
\begin{equation}
\label{maxquo}
\Omega (\rho )\cong \rho \otimes \Theta (\rho ),
\end{equation}
where $\Theta (\rho)$ is a $\widetilde{G}$-module.

Howe duality theorem asserts that $\Theta (\rho)$ is a finitely generated admissible
quasisimple representation of $\widetilde{G}$, and has a unique
irreducible $\widetilde{G}$-quotient, denoted by $\theta (\rho )$:
\[\Theta (\rho)\twoheadrightarrow \theta (\rho).\]

The correspondence
$$\rho \mapsto \theta (\rho )$$
is the Howe quotient correspondence.

\vspace{0.2in}

\noindent {\bf Remark}: If we take $\rho$ to be in the
Casselman-Wallach class \cite[Chapter 11]{Wal}, then we will have
the isomorphism of Casselman-Wallach representations $\Omega (\rho
)\cong \rho \widehat{\otimes}\Theta (\rho )$, where
``$\widehat{\otimes}$'' stands for the completed projective tensor
product.

\section{The coinvariants and the embedding}
\label{Embed}

For the dual pair $(H,G)=({\mathrm{O}(p,q)}, {\mathrm{Sp}(2n,{\mathbb R})})$, let
\begin{align*}
\Omega ^{p,q}= &\ \text{Howe's maximal quotient corresponding to}\\
&\text{ the trivial representation ${1\!\! 1}$ of ${\mathrm{O}(p,q)}$}\\
= &\ {\mathcal S}(M_{p+q,n}({\mathbb R}))_{{\mathrm{O}(p,q)}} \ \ (\text{the space of coinvariants}).
\end{align*}
This is a representation of $\widetilde{G}$.

\vspace{0.1in}

\noindent {\bf Remark}: the continuous dual of $\Omega ^{p,q}$ is ${\mathcal S}^*(M_{p+q,n}({\mathbb R}))^{{\mathrm{O}(p,q)}}$, the space of ${\mathrm{O}(p,q)}$-invariant tempered distributions on $M_{p+q,n}({\mathbb R})$. This is investigated in \cite{KR1} and \cite{Zh}.

\vspace{0.2in}

From now on, we shift notation and will denote $G={\widetilde{\mathrm{Sp}}(2n,{\mathbb R})}$. The case of $n=1$ differs from the general case slightly, but since it is straightforward, we omit it and will assume that $n\geq 2$ throughout this paper. We shall identify
${\widetilde{\mathrm{Sp}}(2n,{\mathbb R})}$ as a set with
\[{\mathrm{Sp}(2n,{\mathbb R})}\times{\mathbb Z}_2=\{(g,\varepsilon):\ g\in{\mathrm{Sp}(2n,{\mathbb R})},\ \varepsilon=\pm 1\}.\]
For $a\in \mathrm{GL}(n,{\mathbb R})$ and $b\in M_{n}({\mathbb R})$ such that $b=b^t$, we let
\begin{eqnarray*}
m_a&=&\left(\begin{array}{cc}
a&0\\
0&(a^{-1})^t
\end{array}\right),\\
&&\\
n_b&=&\left(\begin{array}{cc}
I_n&b\\
0&I_n
\end{array}\right).
\end{eqnarray*}
Let
\[M=\{(m_a,\varepsilon):\ a\in \mathrm{GL}(n,{\mathbb R}),\ \varepsilon=\pm 1\}\]
and
\[N=\{(n_b,1):\ b\in M_n({\mathbb R}),\ b=b^t\}.\]
Then $P=MN$ is a maximal parabolic subgroup of $G$, called the Siegel parabolic.

\vspace{0.2in} Let $\chi:M\longrightarrow{\mathbb C}^\times$ be given by
\[\chi(m_a,\varepsilon)=\varepsilon\cdot\left\{\begin{array}{ll}
i&\mbox{\rm if}\ \det a<0,\\
1&\mbox{\rm if}\ \det a>0.
\end{array}\right.\]
This is a character of $M$ and it is of order $4$. For
$\alpha=0,1,2,3$ and $\sigma\in {\mathbb C}$, let $\chi^\alpha_\sigma$ be the
character of $P$ given by
\begin{equation}
\chi^\alpha_\sigma[(m_a,\varepsilon)(n_b,1)]=|\det a|^\sigma\chi(m_a,\varepsilon)^\alpha.
\end{equation}
Let $I^\alpha(\sigma)$ be the normalized induced representation:
\[I^\alpha(\sigma)=\mbox{\rm Ind}_P^G\chi^\alpha_\sigma .\]
The representation space of ${I^\alpha(\sigma)}$ is
\[\{f\in C^\infty(G):\ f(pg)=\delta^{\frac{1}{2}}(p)
\chi^\alpha_\sigma(p)f(g),\ \forall g\in G, p\in P\},\] and $G$ acts
by right translation:
\[g\cdot f(h)=f(hg),\ \ \ \ (g,h\in G).\]
Here $\delta$ denotes the modular function of $P$, and is given by
\[\delta[(m_a,\varepsilon)(n_b,1)]=|\det a|^{2\rho_n}\]
and
\[\rho_n=\frac{n+1}{2}.\]
When $\alpha=0$ or $2$, the representation ${I^\alpha(\sigma)}$ descends to a representation of the
linear group ${\mathrm{Sp}(2n,{\mathbb R})}$.

\vspace{0.2in}

Define the map $\psi _{p,q}: {\mathcal S}(M_{p+q,n}({\mathbb R}))\mapsto C^\infty(G)$ by
\[\psi _{p,q}(f)(g)=(\omega (g)f)(0), \ \ f\in {\mathcal S}(M_{p+q,n}({\mathbb R})), \ g\in G.\]
From the well-known formula of the oscillator representation $\omega$, we see that
\[\psi _{p,q}: {\mathcal S}(M_{p+q,n}({\mathbb R}))\mapsto I^\alpha(\sigma),\]
where
\begin{equation}
\label{br}
\sigma = \frac{p+q}{2}-\rho_n, \ \ \text{and} \ \ \alpha \equiv p-q \ (\mathrm{mod}\  4).
\end{equation}

The following is the fundamental result of Kudla-Rallis.

\vspace{0.2in}

\begin{thm} {\rm (\cite{KR1})} The map $\psi _{p,q}$ induces a topological embedding with closed image:
\[\Omega ^{p,q}\hookrightarrow I^\alpha(\sigma)\]
\end{thm}

\vspace{0.2in}

\noindent {\bf Remark}:  Analogous results hold for other classical groups. See \cite{Zh}.

\section{Degenerate principal series: the transition coefficients}
\label{TC}

Fix a maximal compact subgroup $K_1\simeq {\rm U}(n)$ of ${\mathrm{Sp}(2n,{\mathbb R})}$ and let $K$ be the inverse image of $K_1$ in ${\widetilde{\mathrm{Sp}}(2n,{\mathbb R})}$, thus a maximal compact subgroup of ${\widetilde{\mathrm{Sp}}(2n,{\mathbb R})}$. The Lie algebra of ${\mathrm{Sp}(2n,{\mathbb R})}$ has a Cartan decomposition
\[\mathfrak{sp}(2n,{\mathbb R})={\mathfrak k} \oplus {\mathfrak p}\]
where ${\mathfrak k}$ is the Lie algebra of $K_1$ and is isomorphic to ${\mathfrak u}(n)$.

Let
\[{\Lambda^+_n}=\{(\lambda_1,...,\lambda_n)\in{\mathbb Z}^n:
\ \lambda_1\geq\lambda_2\geq\cdot\cdot\cdot\geq \lambda_n\},\]
and
\[{\mbox{\rm\bf 1}}=(1,...,1)\in{\Lambda^+_n}.\]
Also let $e_j=(\overbrace{0,...,0,1}^{j},0,...,0)$, for $1\leq j\leq n$.

We have the $K$-type decomposition:
\begin{equation}
{I^\alpha(\sigma)}|_K=\oplus_{\lambda\in{\Lambda^+_n}}V_{2\lambda+{\textstyle \frac{\alpha}{2}}{\mbox{\rm\bf 1}}},
\end{equation}
where $V_{2\lambda+\frac{\alpha}{2}{\mbox{\rm\bf 1}}}$ is an
irreducible $K$-module with highest weight $2\lambda+{\textstyle
\frac{\alpha}{2}}{\mbox{\rm\bf 1}}$, for $\lambda\in{\Lambda^+_n}$.

For $\mu=2\lambda+\frac{\alpha}{2}{\mbox{\rm\bf 1}}$, let $V_{\mu}$
be a $K$-type in ${I^\alpha(\sigma)}$ and fix $\gamma_{\mu}$ to the
unique (up to a multiple) $K$-highest weight vector in $V_{\mu}$. We
consider the tensor product ${\mathfrak p}_{\mathbb C}\otimes
V_{\mu}$ of ${\mathfrak k}_{\mathbb C}$ modules, and the ${\mathfrak
k}_{\mathbb C}$ map
\[m:{\mathfrak p}_{\mathbb C}\otimes V_{\mu} \longrightarrow {I^\alpha(\sigma)}_K\]
\[m(p\otimes v)=p\cdot v.\]
For each $1\leq j\leq n$, there exists an element $X_j$ in
${\mathcal U} ({\mathfrak {sp}}(2n,{\mathbb C}))$ (the universal
enveloping algebra of the complexified Lie algebra) with the
property that $X_j\cdot \gamma_\mu$ is the image of the unique
${\mathfrak k}_{\mathbb C}$ highest weight vector of weight
$\mu+2e_j$ in ${\mathfrak p}_{\mathbb C}\otimes V_{\mu}$, valid for
all $\mu$. C.f. \cite[Section 3]{L2}. Thus $X_j\cdot \gamma_\mu $ is
a multiple of the unique ${\mathfrak k}_{\mathbb C}$ highest weight
vector $\gamma_{\mu+2e_j}$ in $V_{\mu+2e_j}$. We call this multiple
the transition coefficient from $V_{\mu}$ to $V_{\mu +2e_j}$.
Likewise there is a transition coefficient from $V_{\mu}$ to $V_{\mu
-2e_j}$. Note that the transition coefficients depend on the choice
of $\gamma_{\mu}$'s.

The main use of the transition coefficients is as follows: if the
transition coefficient from $V_{\mu}$ to $V_{\mu +2e_j}$ is nonzero,
then $V_{\mu +2e_j}$ is in the submodule generated by $V_{\mu}$, and
vice versa. We shall indicate by the symbol: $V_\mu\rightarrow V_{\mu+2e_j}$; 
If this transition coefficient is zero, then $\sum_{\mu '_j\leq \mu
_j}V_{\mu '}$ is a submodule of ${I^\alpha(\sigma)}$ and we say that
we have a barrier to block $j$th rightward movement
$V_\mu\rightarrow V_{\mu+2e_j}$. Likewise for the transition
coefficient from $V_{\mu}$ to $V_{\mu -2e_j}$. To summarize, the
collection of transition coefficients allows us to determine the
lattice of submodules in ${I^\alpha(\sigma)}$, explicitly.

By constructing an appropriate choice of ${\mathfrak k}_{\mathbb C}$
highest weight vector $\gamma_\mu$ for all $\mu$ (c.f.
\cite[Section3]{L2}), the transition coefficients $A_j^+(\mu)$ from
$V_{\mu}$ to $V_{\mu+2e_j}$, and the transition coefficients
$A_j^-(\mu)$ from $V_{\mu}$ to $V_{\mu-2e_j}$, are computed to be
\begin{equation}
\label{tc1}
A_j^+(\mu)=B^+_j-2\lambda _j, \ \ \ A_j^-(\mu)=2\lambda _j-B^-_j,
\end{equation}
where
\begin{equation}
\label{tc2}
B^+_j= -\sigma-\rho_n-{\textstyle \frac{\alpha}{2}}+j-1, \ \ \ B^-_j= \sigma-\rho_n-{\textstyle \frac{\alpha}{2}}+j+1.
\end{equation}

\vspace{0.2in}
\begin{cor}\label{red} {\rm (\cite{KR1})}
${I^\alpha(\sigma)}$ is irreducible if and only if
\[\sigma+\rho_n+{\textstyle \frac{\alpha}{2}}\not\in{\mathbb Z}.\]
\end{cor}

\vspace{0.2in} The explicit form of transition coefficients also allows us to determine the complementary series. We will
follow the method used in \cite{L2}.
Each $K$-type
$V_{2\lambda+\frac{\alpha}{2}{\mbox{\rm\bf 1}}}$ of ${I^\alpha(\sigma)}$ has a $K$-invariant inner product given by
\[\langle f_1,f_2\rangle_\lambda=\int_Kf_1(k)\overline{f_2(k)}dk.\]
Since $V_{2\lambda+\frac{\alpha}{2}{\mbox{\rm\bf 1}}}$ is an irreducible $K$-module, any $K$-invariant inner product on
$V_{2\lambda+\frac{\alpha}{2}{\mbox{\rm\bf 1}}}$ is a multiple of $\langle .,.
\rangle_\lambda$. Thus if $\langle .,.\rangle$ is a ${\widetilde{\mathrm{Sp}}(2n,{\mathbb R})}$-invariant
inner product on ${I^\alpha(\sigma)}$ then there exists positive constants
$\{c_\lambda\}_{\lambda\in{\Lambda^+_n}}$  such that
\[\langle f_1,f_2\rangle=
c_\lambda \langle f_1,f_2\rangle_\lambda, \hspace{0.5in}\forall
f_1,f_2\in V_{2\lambda+{\textstyle \frac{\alpha}{2}}{\mbox{\rm\bf 1}}}.\] Since
the $K$-types of ${I^\alpha(\sigma)}$ are mutually orthogonal with respect
$\langle.,.\rangle$, $\langle.,.\rangle$ is completely determined by
the constants $\{c_\lambda\}$. Using similar arguments as the
$U(n,n)$ case (see \cite[section 9]{L1}), we obtain the following:

\begin{lemma} The inner product on ${I^\alpha(\sigma)}$ defined by the constants
$\{c_\lambda\}_{\lambda\in{\Lambda^+_n}}$ is ${\widetilde{\mathrm{Sp}}(2n,{\mathbb R})}$-invariant if and only if
\[(-\sigma-\rho_n-{\textstyle \frac{\alpha}{2}}-2\lambda_j+j-1)c_{\lambda+e_j}
+(-\overline{\sigma}+\rho_n+{\textstyle \frac{\alpha}{2}}+
2\lambda_j-j+1) c_\lambda=0\] for all $\lambda\in{\Lambda^+_n}$ and all
$1\leq j\leq n$.
\end{lemma}

\vspace{0.2in} We let
\begin{equation*}
N_{\lambda,j}=\frac{-\sigma-\rho_n-{\textstyle
\frac{\alpha}{2}}-2\lambda_j+j-1}
{-\overline{\sigma}+\rho_n+{\textstyle \frac{\alpha}{2}}+
2\lambda_j-j+1} =-\frac{c_\lambda}{c_{\lambda+e_j}}.\label{nlam}
\end{equation*}
Then ${I^\alpha(\sigma)}$ is unitarizable if and only if $N_{\lambda,j}
<0$ for all $\lambda\in{\Lambda^+_n}$ and for all $j$.

\vspace{0.2in}
Note that
\[N_{\lambda,j}=\frac{-\sigma-(\rho_n+{\textstyle \frac{\alpha}{2}}+2\lambda_j-j+1)}{
-\overline{\sigma}+(\rho_n+{\textstyle \frac{\alpha}{2}}
+2\lambda_j-j+1)}.\] We write $\xi=\rho_n+{\textstyle
\frac{\alpha}{2}}+2\lambda_j-j+1$. Then
$N_{\lambda,j}=(-\sigma-\xi)/(-\overline{\sigma}+\xi)$. Thus
$N_{\lambda,j}$ is real for all $\lambda$ and for $j$ if and only if
either $\mbox{Re}(\sigma) =0$ or $\sigma$ is real. The case
$\mbox{Re}(\sigma)=0$ corresponds to the unitary axis. If $\sigma$
is real, then
\[N_{\lambda,j}=\frac{-\sigma-\xi}{
-\overline{\sigma}+\xi}<0\Longleftrightarrow |\sigma|<|\xi|.\]

The minimum value of $|\xi|$ is ${\textstyle \frac{1}{2}}$ if $n+\alpha $ is even, and $0$ if $n+\alpha $ is odd. This leads to the following

\begin{thm} \label{coms} If $n+\alpha$ is even,
then ${I^\alpha(\sigma)}$ is unitarizable for $|\sigma|<{\textstyle \frac{1}{2}}$.
\end{thm}

\section{Subquotients of $I^\alpha(\sigma)$}
\label{Subquo}

In this section, we shall give a detailed description of the module
structure of ${I^\alpha(\sigma)}$ when it is reducible. We shall
describe all the irreducible subquotients of ${I^\alpha(\sigma)}$
and determine which of them are unitarizable, i.e., possess a
$G$-invariant positive
definite inner product. We also describe the socle series and module
diagram of ${I^\alpha(\sigma)}$.

\vspace{0.2in} Let
\begin{equation}
\label{dtsigma}
\tilde{\sigma} =\sigma+\rho_n+{  \frac{\alpha}{2}} =\sigma+{
\frac{n+1+\alpha}{2}}.
\end{equation}
By Corollary \ref{red},
$I^\alpha(\sigma)$ is irreducible if and only if $\tilde{\sigma}\not\in{\mathbb Z}$. Thus
we shall assume that $\tilde{\sigma} \in{\mathbb Z}$ throughout this section.

Let $\mu=2\lambda+(\alpha/2){\mbox{\rm\bf 1}}_n$, as before. Then by the
results on transition coefficients in Section \ref{TC}, we have
\[V_\mu\rightarrow V_{\mu+2e_j}\Longleftrightarrow 2\lambda_j\neq B^+_j\hspace{0.25in}
\mbox{and}\hspace{0.25in} V_\mu\rightarrow V_{\mu-2e_j}\Longleftrightarrow
2\lambda_j\neq B^-_j\] where
\[B^+_j=B^+_j(\alpha,\sigma)=-\sigma-\frac{n+\alpha+1}{2}+j-1=-\tilde{\sigma} +j-1\]
and
\[B^-_j=B^-_j(\alpha,\sigma)=\sigma-\frac{n+\alpha+1}{2}+j+1=\tilde{\sigma}-(n+\alpha)+j.\]

As discussed in Section \ref{TC}, we may determine the lattice of
submodules of $I^\alpha(\sigma)$ from the transition coefficients
$\{B^{\pm}_j\}_{1\leq j\leq n}$. Schematically it goes as follows:
\begin{enumerate}
\item Represent $V_\mu$ by the point $2\lambda\in{\mathbb R}^n$ with
standard coordinates $(x_1,...,x_n)$.
\item Barrier to block $j$-th rightward movement $V_\mu\rightarrow V_{\mu+2e_j}$
is at hyperplane $\ell^+_j: x_j=B^+_j$.
\item Barrier to block $j$-th leftward movement $V_\mu\rightarrow V_{\mu-2e_j}$
is at hyperplane $\ell^-_j: x_j=B^-_j$.
\item Analyze barriers systematically (Subsections \ref{Subodd} and \ref{Subeven}). 
\end{enumerate}
Note that the barrier $\ell^+_j: x_j=B^+_j$ is effective if and only
if it cuts at an ``even" point, i.e. $B^+_j\in 2{\mathbb Z}$. The
parity of $B^+_j$ depends on the parity of $\tilde{\sigma}
=\sigma+\frac{n+\alpha+1}{2}$ and the parity of $j$. Similarly for
the  barrier $\ell^-_j $, but the parity of $B^-_j$ also depends on
the parity of $n+\alpha$. So our analysis will be divided into 4
cases.

\vspace{0.2in}
\begin{center}
\begin{tabular}{|c|c|c|}\hline
  $\tilde{\sigma} \backslash n+\alpha $& odd & even \\ \hline odd & Case 1a & Case 2a\\
  \hline even &
Case 1b   & Case 2b\\  \hline
\end{tabular}
\end{center}

\vspace{0.2in} \noindent We also define the ``gap" between $\ell^+_j$ and $\ell^-_j$  by
\begin{equation}\label{gap}
\mbox{gap}=B^+_j-B^-_j=-2\sigma-2.
\end{equation}

From the lattice of submodules of $I^\alpha(\sigma)$, we may give
its module structure a diagrammatic representation, called the
module diagram. It is a directed simple graph $\mathcal{G}$, defined
as follows: The vertex set of $\mathcal{G}$ is the set of all
irreducible subquotients in $I^\alpha(\sigma)$, to be identified
with its collection of $K$-types. There is a directed edge from the
node $R_1$ to the node $R_2$ if and only if there are submodules $U$
and $V$ of $I^\alpha(\sigma)$ such that $V\subseteq U$ and there is
a nonsplit exact sequence of infinitesimal
$\widetilde{\mathrm{Sp}}(2n,{\mathbb R})$-modules $0\rightarrow
R_2\rightarrow U/V\rightarrow R_1\rightarrow 0$. We shall also
arrange the nodes in $\mathcal{G}$ in such a way that all the edges
are directed downward. Then one can recover the lattice of
submodules of  $I^\alpha(\sigma)$ from the graph $\mathcal{G}$. For
a more detailed explanation on module diagram, see Section 7 of
\cite{L1}.

\vspace{0.2in}

\subsection{Subquotients of ${I^\alpha(\sigma)}$: $n+\alpha $ odd}\label{Subodd}\[\]

\noindent{\bf Case 1: $n+\alpha$ odd.}

\vspace{0.2in}
 Since $\tilde{\sigma} \in{\mathbb Z}$ and $\frac{n+\alpha+1}{2}\in{\mathbb Z}$, we have $\sigma\in{\mathbb Z}$. It follows from this and equation \eqref{gap} that the gap $B^+_j-B^-_j$ is even, i.e.
 $B^+_j\equiv B^-_j(\mathrm{mod}\ 2)$. This means that for each $j$, either both barriers
$\ell^+_j$ and $\ell^-_j$ are effective, or both are not effective.

\vspace{0.2in} \noindent We consider two subcases:
\[\mbox{{\em Case 1a:}  $\tilde{\sigma} =\sigma+\frac{n+\alpha+1}{2}$ is odd\hspace{0.5in} {\em Case 1b:} $\tilde{\sigma} =\sigma+\frac{n+\alpha+1}{2}$ is even.}\]

\vspace{0.2in}\noindent For Case 1a, We have
\[B^+_j= -(\sigma+\frac{n+\alpha+1}{2})+j-1\in 2{\mathbb Z}\Longleftrightarrow \mbox{$j$ is even.}\]
Hence the effective barriers are  $\ell^+_{2i}$ and $\ell^-_{2i}$
for all $1\leq i\leq n_0/2$, where
\begin{equation}
\label{dn0}
n_0=2\left[\frac{n}{2}\right]=\mbox{largest even integer less than or equal to $n$.}
\end{equation}
%Also,
%\[
%B^-_j= \sigma-\frac{n+\alpha+1}{2}+j+1=(\sigma+\frac{n+\alpha+1}{2})
%-(n+\alpha+1)+j+1\in 2{\mathbb Z}\Longleftrightarrow \mbox{$j$ is even}\]

\vspace{0.2in}\noindent But for Case 1b, We have
\[B^+_j= -(\sigma+\frac{n+\alpha+1}{2})+j-1\in 2{\mathbb Z}\Longleftrightarrow \mbox{$j$ is odd.}\]
So the effective barriers are  $\ell^+_{2i-1}$ and $\ell^-_{2i-1}$
for all $1\leq i\leq (n_1+1)/2$, where
\begin{equation}
\label{dn1}
n_1=2\left[\frac{n+1}{2}\right]-1=\mbox{largest odd integer less than or equal to $n$.}
\end{equation}

\vspace{0.2in} \noindent {\bf Case 1a:}  $\tilde{\sigma} =\sigma+\frac{n+\alpha+1}{2}$ is odd

\vspace{0.2in} \noindent Let  $i$ and $j$ be such that
 $0\leq i+j\leq n_0/2$. We will define  $R_{ij}(n,\sigma,\alpha)$ as
 follows:

 \vspace{0.2in}\begin{enumerate}
 \item[(i)] For $\sigma\leq -1$,
 let $R_{ij}(n,\sigma,\alpha)$ be the set of $2\lambda$ such that
\begin{equation}
2\lambda_{2i}\geq B^+_{2i+2}\geq 2\lambda_{2i+2}
\end{equation}
\begin{equation}
2\lambda_{n_0-2j}\geq B^-_{n_0-2j} \geq 2\lambda_{n_0-2j+2}
\end{equation}
where
\[
B^+_{2i+2}=-\sigma-\frac{n+\alpha+1}{2}+2i+1=-\sigma-\frac{n+\alpha-1}{2}+2i,
\]
\[
B^-_{n_0-2j}=\sigma-\frac{n+\alpha+1}{2}+n_0-2j+1=\sigma-\frac{n+\alpha-1}{2}+n_0-2j.\]

\vspace{0.2in} \item[(ii)] For $\sigma\geq 0$,
 let $R_{ij}(n,\sigma,\alpha)$ be the set of $2\lambda$ such that
\begin{equation}\label{bcon1}
2\lambda_{2i}\geq B^-_{2i}\geq 2\lambda_{2i+2}
\end{equation}
\begin{equation}\label{bcon2}
2\lambda_{n_0-2j}\geq B^+_{n_0-2j+2} \geq 2\lambda_{n_0-2j+2}
\end{equation}
 where
\[B^-_{2i}=\sigma-\frac{n+\alpha+1}{2}+2i+1=\sigma-\frac{n+\alpha-1}{2}+2i,\]
\[B^+_{n_0-2j+2}=-\sigma-\frac{n+\alpha+1}{2}+n_0-2j+1=-\sigma-\frac{n+\alpha-1}{2}+n_0-2j.\]
\end{enumerate}
Note that the set  $R_{ij}(n,\sigma,\alpha)$ defined
above may be empty. When it is nonempty, it will be identified with
  the direct sum of all the $K$-representations
$V_{2\lambda+\frac{\alpha}{2}{\mbox{\rm\bf 1}}}$ in
$I^\alpha(\sigma)$ such that $2\lambda\in  R_{ij}(n,\sigma,\alpha)$.

\vspace{0.2in}\noindent{\bf Remark:} When $\sigma=0$ and
$i+j=n_0/2$, the two conditions \eqref{bcon1} and \eqref{bcon2}
coincide. More precisely,
\[R_{ij}(n,0,\alpha)=\{2\lambda:\ 2\lambda_{2i}\geq B^-_{2i}\geq
2\lambda_{2i+2}\}.\]

\vspace{0.2in}

In what follows, $\oplus$ denotes sum in the Grothendieck group of
$G$-modules. If we have a direct sum of submodules, we will state it
explicitly.

\vspace{0.2in}\begin{thm} {\rm (Case 1a)} Assume that $n+\alpha$ is odd, $\sigma$ is an integer and
  $\tilde{\sigma} =\sigma+\frac{n+\alpha+1}{2}$ is odd.

\begin{enumerate}
\item[(a)] If $\sigma\leq -1$, then
\[{I^\alpha(\sigma)}=\bigoplus\left\{R_{ij}(n,\sigma,\alpha):\ r_1\leq i+j\leq
\frac{n_0}{2}\right\},\]  where
\[r_1=\max\left(\frac{n_0}{2}+\sigma,0\right).\]
In this case, the module diagram of ${I^\alpha(\sigma)}$ can be obtained
from Figure 1 by removing those $R_{ab}(n,\sigma,\alpha)$ which are
empty. In particular, the socle series of ${I^\alpha(\sigma)}$ is given by
\[\mbox{\rm Soc}^l({I^\alpha(\sigma)})=\left\{\begin{array}{ll}
\bigoplus_{ r_1\leq i+j\leq r_1+l-1 } R_{ij}(n,\sigma,\alpha) &1\leq l\leq
\frac{n_0}{2}-r_1,\\
&\\ {I^\alpha(\sigma)}&l\geq \frac{n_0}{2} -r_1+1.
\end{array}\right.\]
An irreducible constituent $R_{ij}(n, \sigma, \alpha)$ of ${I^\alpha(\sigma)}$ is
unitarizable if and only if
$-\frac{n_0}{2}\leq \sigma\leq -1  $ and $i+j=r_1$.

\vspace{0.2in}\item[(b)] If $\sigma=0$, then
\[I^\alpha(0)=\bigoplus_{i+j=n_0/2}R_{ij}(n,0,\alpha)\]
is a direct sum of irreducible unitary submodules.

\vspace{0.2in}
\item[(c)] If $\sigma\geq 1$,
then
\[{I^\alpha(\sigma)}=\bigoplus\left\{R_{ij}(n,\sigma,\alpha):\ r_2\leq i+j\leq
\frac{n_0}{2}\right\},\]  where
\[r_2=\max\left(\frac{n_0}{2}-\sigma,0\right).\]
In this case,
the module diagram of ${I^\alpha(\sigma)}$ can be obtained
from Figure 2 by removing those $R_{ab}(n,\sigma,\alpha)$ which are
empty. In particular, the socle series of ${I^\alpha(\sigma)}$ is given by
\[\mbox{\rm Soc}^l({I^\alpha(\sigma)})=\left\{\begin{array}{ll}
\bigoplus_{ \frac{n_0}{2}-l+1  \leq i+j\leq \frac{n_0}{2}    }R_{ij}(n,\sigma,\alpha)&1\leq l\leq
\frac{n_0}{2}-r_2,\\
&\\ {I^\alpha(\sigma)}&l\geq \frac{n_0}{2} -r_2+1.
\end{array}\right.\]
An irreducible constituent $R_{ij}(n, \sigma, \alpha)$ of ${I^\alpha(\sigma)}$ is
unitarizable if and only if
\   $1\leq \sigma\leq
\frac{n_0}{2}  $ and $i+j=r_2$.

\end{enumerate}
\end{thm}

\vspace{0.2in} \noindent {\bf Case 1b:}  $\tilde{\sigma} =\sigma+\frac{n+\alpha+1}{2}$ is even

\vspace{0.2in} Let  $i$ and $j$ be such that
 $0\leq i+j\leq (n_1+1)/2$. We will define  $R_{ij}(n,\sigma,\alpha)$ as
 follows:

 \vspace{0.2in}\begin{enumerate}
 \item[(i)] For $\sigma\leq -1$,
 let $R_{ij}(n,\sigma,\alpha)$ be the set of $2\lambda$ such that
 \begin{equation}\label{bcon11}
2\lambda_{2i-1}\geq B^+_{2i+1}\geq 2\lambda_{2i+1}
\end{equation}
\begin{equation}\label{bcon22}
2\lambda_{n_1-2j}\geq B^-_{n_1-2j} \geq
2\lambda_{n_1-2q+2}\end{equation} where
\[B^+_{2i+1}=-\sigma-\frac{n+\alpha+1}{2}+2i,\]
\[B^-_{n_1-2j}=\sigma-\frac{n+\alpha+1}{2}+n_1-2j+1=
\sigma-\frac{n+\alpha-1}{2}+n_1-2j.\]

\vspace{0.2in} \item[(ii)] For $\sigma\geq 0$,
 let $R_{ij}(n,\sigma,\alpha)$ be the set of $2\lambda$ such that
 \begin{equation}
2\lambda_{2i-1}\geq B^-_{2i-1}\geq 2\lambda_{2i+1}
\end{equation}
\begin{equation}
2\lambda_{n_1-2j}\geq B^+_{n_1-2j+2} \geq
2\lambda_{n_1-2j+2}
\end{equation} where
\[B^-_{2i-1}=\sigma-\frac{n+\alpha+1}{2}+2i,\]
\[B^+_{n_1-2j+2}=-\sigma-\frac{n+\alpha+1}{2}+n_1-2j+1
=-\sigma-\frac{n+\alpha-1}{2}+n_1-2j.\]
\end{enumerate}

\vspace{0.2in}\noindent{\bf Remark:}   When $\sigma=0$ and $i+j=(n_1+1)/2$, the two
conditions \eqref{bcon11} and \eqref{bcon22} coincide, that is, $
2\lambda_{n_1-2j}=2\lambda_{(2i+2j-1)-2j}=2\lambda_{2i-1}$,
\[B^-_{n_1-2j+2}=-\frac{n+\alpha-1}{2}+(2i+2j-1)-2j
=-\frac{n+\alpha+1}{2}+2i=B^-_{2i-1},\] and $
2\lambda_{n_1-2j+2}=2\lambda_{(2i+2j-1)-2j+2}=2\lambda_{2j+1}$.
Precisely,
\[R_{ij}(n,0,\alpha)=\{2\lambda:\ 2\lambda_{2i-1}\geq B^-_{2i-1}\geq
2\lambda_{2i+1}\}.\]

\vspace{0.2in}\begin{thm} {\rm (Case 1b)} Assume that $n+\alpha$ is odd, $\sigma$ is an integer and
  $\tilde{\sigma} =\sigma+\frac{n+\alpha+1}{2}$ is even.

\begin{enumerate}
\item[(a)] If $\sigma\leq -1$, then
\[{I^\alpha(\sigma)}=\bigoplus\left\{R_{ij}(n,\sigma,\alpha):\ r_1\leq i+j\leq
\frac{n_1+1}{2}\right\},\]  where
\[r_1=\max\left(\frac{n_1+1}{2}+\sigma,0\right).\]
In this case, the module diagram of ${I^\alpha(\sigma)}$ can be obtained
from Figure 1 by removing those $R_{ab}(n,\sigma,\alpha)$ which are
empty. In particular, the socle series of ${I^\alpha(\sigma)}$ is given by
\[\mbox{\rm Soc}^l({I^\alpha(\sigma)})=\left\{\begin{array}{ll}
\bigoplus_{ r_1\leq i+j\leq r_1+l-1 } R_{ij}(n,\sigma,\alpha) &1\leq l\leq
\frac{n_1+1}{2}-r_1,\\
&\\ {I^\alpha(\sigma)}&l\geq \frac{n_1+1}{2} -r_1+1.
\end{array}\right.\]
An irreducible constituent $R_{ij}(n, \sigma, \alpha)$ of ${I^\alpha(\sigma)}$ is
unitarizable if and only if
\begin{enumerate}
\item[(i)]   $
-\frac{n_1+1}{2}\leq \sigma\leq -1  $ and $i+j=r_1$; or
\item[(ii)] $n$ is odd, $\alpha\in\{0,2\}$ and
$(i,j)\in\left\{\left(\frac{n+1}{2},0\right),\left(0,\frac{n+1}{2}\right)\right\}$.
\end{enumerate}

\vspace{0.2in}\item[(b)] If $\sigma=0$, then
\[I^\alpha(0)=\bigoplus_{i+j=\frac{n_1+1}{2}}R_{ij}(n,0,\alpha)\]
is a direct sum of irreducible unitary submodules.

\vspace{0.2in}
\item[(c)] If $\sigma\geq 1$,
then
\[{I^\alpha(\sigma)}=\bigoplus\left\{R_{ij}(n,\sigma,\alpha):\ r_2\leq i+j\leq
\frac{n_1+1}{2}\right\},\]  where
\[r_2=\max\left(\frac{n_1+1}{2}-\sigma,0\right).\]
In this case,
the module diagram of ${I^\alpha(\sigma)}$ can be obtained
from Figure 2 by removing those $R_{ab}(n,\sigma,\alpha)$ which are
empty. In particular, the socle series of ${I^\alpha(\sigma)}$ is given by
\[\mbox{\rm Soc}^l({I^\alpha(\sigma)})=\left\{\begin{array}{ll}
\bigoplus_{ \frac{n_1+1}{2}-l+1  \leq i+j\leq \frac{n_1+1}{2}    }R_{ij}(n,\sigma,\alpha)&1\leq l\leq
\frac{n_1+1}{2}-r_2,\\
&\\ {I^\alpha(\sigma)}&l\geq \frac{n_1+1}{2} -r_2+1.
\end{array}\right.\]
An irreducible constituent $R_{ij}(n,\sigma,\alpha)$ of ${I^\alpha(\sigma)}$ is
unitarizable if and only if
\begin{enumerate}
\item[(i)]   $1\leq \sigma\leq
\frac{n_1+1}{2}  $ and $i+j=r_2$; or
\item[(ii)] $n$ is odd, $\alpha\in\{0,2\}$   and
$(i,j)\in\left\{\left(\frac{n+1}{2},0\right),\left(0,\frac{n+1}{2}\right)\right\}$.
\end{enumerate}

\end{enumerate}
\end{thm}

\vspace{0.2in} The module diagram of ${I^\alpha(\sigma)}$ for $\sigma\leq -1$ and $\sigma\geq 1$
 are given in Figure 1 and Figure 2 below.
Here
\begin{equation}
\label{defk}
k=\left\{\begin{array}{ll}\displaystyle
\frac{n_0}{2}=[\frac{n}{2}]&\mbox{\hspace{0.5in} if $\tilde{\sigma}$ is odd,}\\
&\\  \displaystyle\frac{n_1+1}{2}=[\frac{n+1}{2}]&\mbox{\hspace{0.5in} if $\tilde{\sigma}$ is even.}
\end{array}\right.
\end{equation}

and we write a constituent  $R_{ij}(n,\sigma,\alpha)$ simply as
$R_{ij}$.

\vspace{0.2in} {\tiny \setlength{\unitlength}{0.8mm}

\begin{center}
\begin{picture}(120,65)(0,5)

\put(0,60){\makebox(10,10){$\scriptstyle
R_{k,0} $}}
\put(20,60){\makebox(10,10){$\scriptstyle
R_{k-1,1}   $}}

\put(40,65){$\cdot$} \put(50,65){$\cdot$} \put(60,65){$\cdot$}
\put(70,65){$\cdot$} \put(80,65){$\cdot$}

\put(90,60){\makebox(10,10){$\scriptstyle R_{1,k-1}  $}}
\put(110,60){\makebox(10,10){$\scriptstyle R_{0,k}  $}}

\put(25,35){\makebox(10,10){$\scriptstyle R_{2,0} $}}
\put(55,35){\makebox(10,10){$\scriptstyle R_{1,1} $}}
\put(85,35){\makebox(10,10){$\scriptstyle R_{0,2} $}}

\put(27,44){\line(-1,1){5}} \put(33,44){\line(1,1){5}}
\put(57,44){\line(-1,1){5}} \put(63,44){\line(1,1){5}}
\put(87,44){\line(-1,1){5}} \put(93,44){\line(1,1){5}}

\put(40,20){\makebox(10,10){$\scriptstyle R_{1,0} $}} \put(70,20){\makebox(10,10){$\scriptstyle R_{0,1}  $}}

\put(42,29){\line(-1,1){7}} \put(48,29){\line(1,1){7}}

\put(72,29){\line(-1,1){7}} \put(78,29){\line(1,1){7}}

\put(20,51){.} \put(18,53){.} \put(16,55){.}

\put(14,57){\line(-1,1){5}} \put(16,57){\line(1,1){5}}
\put(30,57){\line(-1,1){5}}

 \put(100,51){.} \put(102,53){.} \put(104,55){.}
\put(106,57){\line(1,1){5}} \put(104,57){\line(-1,1){5}}
\put(90,57){\line(1,1){5}}

\put(57,14){\line(-1,1){7}} \put(63,14){\line(1,1){7}}
\put(55,5){\makebox(10,10){$\scriptstyle R_{0,0} $}}
\end{picture}

Figure 1: Module diagram for ${I^\alpha(\sigma)}$ $(\sigma\leq -1)$
\end{center}

{\vspace{0.05in}\noindent}

\vspace{0.2in} \setlength{\unitlength}{0.8mm}

\begin{center}
\begin{picture}(120,65)(0,5)

\put(0,10){\makebox(10,10){$\scriptstyle R_{k,0}   $}}
\put(20,10){\makebox(10,10){$\scriptstyle R_{k-1,1}$}}

\put(40,15){$\cdot$} \put(50,15){$\cdot$} \put(60,15){$\cdot$}
\put(70,15){$\cdot$} \put(80,15){$\cdot$}

\put(90,10){\makebox(10,10){$\scriptstyle R_{1,k-1}$}}
\put(110,10){\makebox(10,10){$\scriptstyle R_{0,k}$}}

%\put(13,23){$\cdot$} \put(17,27){$\cdot$} \put(21,31){$\cdot$}
%\put(99,31){$\cdot$} \put(103,27){$\cdot$} \put(107,23){$\cdot$}

\put(25,35){\makebox(10,10){$\scriptstyle R_{2,0}$}}
\put(55,35){\makebox(10,10){$\scriptstyle R_{1,1}$}}
\put(85,35){\makebox(10,10){$\scriptstyle R_{0,2}$}}

\put(40,50){\makebox(10,10){$\scriptstyle R_{1,0}$}}
\put(70,50){\makebox(10,10){$\scriptstyle R_{0,1}$}}

\put(47,60){\line(1,1){7}}
\put(69,60){\line(-1,1){7}}

\put(31,44){\line(1,1){7}} \put(55,44){\line(-1,1){7}}
\put(61,44){\line(1,1){7}} \put(84,44){\line(-1,1){7}}

\put(24,37){\line(-1,-1){5}} \put(29,37){\line(1,-1){5}}
\put(55,37){\line(-1,-1){5}} \put(60,37){\line(1,-1){5}}
\put(85,37){\line(-1,-1){5}} \put(91,37){\line(1,-1){5}}

\put(13,26){.} \put(15,28){.} \put(17,30){.}
\put(102,26){.} \put(100,28){.} \put(98,30){.}

\put(6,19){\line(1,1){5}} \put(109,19){\line(-1,1){5}}
\put(23,19){\line(-1,1){5}} \put(25,19){\line(1,1){5}}
\put(92,19){\line(-1,1){5}} \put(94,19){\line(1,1){5}}

\put(55,65){\makebox(10,10){$\scriptstyle R_{0,0}$}}

\end{picture}

Figure 2: Module diagram for ${I^\alpha(\sigma)}$ $(\sigma\geq 1)$

\end{center}}

\subsection{Subquotients of ${I^\alpha(\sigma)}$: $n+\alpha $ even}
\label{Subeven}
\[
\]

\vspace{0.2in}\noindent{\bf Case 2: $n+\alpha$ even.}

\vspace{0.2in} Recall that $\tilde{\sigma} =\sigma+\frac{n+\alpha+1}{2}\in{\mathbb Z}$. So in this case, $\sigma+{\textstyle \frac{1}{2}} \in{\mathbb Z}$, and consequently the
gap between the barriers $\ell^+_j$ and $\ell^-_j$ given by
\[\mbox{gap}=-2\sigma-2\]
is odd. It follows that along each coordinate axis, exactly one barrier is effective. More precisely, for each $1\leq j\leq n$, either $\ell^+_j$
is effective or $\ell^-_j$ is effective, but not both. Again, we
consider two subcases:
\[\mbox{{Case 2a:}  $\tilde{\sigma} =\sigma+\frac{n+\alpha+1}{2}$ is odd\hspace{0.5in} { Case 2b:}
$\tilde{\sigma} = \sigma+\frac{n+\alpha+1}{2}$ is even.}\] For Case 2a, $\ell^+_j$ is
effective when $j$ is even, and $\ell^-_j$ is effective when $j$ is
odd. Case 2b is the other way round.

\vspace{0.2in}\noindent{\bf Case 2a:}  $\tilde{\sigma} =\sigma+\frac{n+\alpha+1}{2}$ is odd

Since exactly one barrier is effective along each coordinate axis, the barrier partitions the $K$-types into $2$ subsets.
For $1\leq r\leq \frac{n_1+1}{2}+1$, let
\begin{eqnarray*}
X^r_1&=&\{2\lambda: 2\lambda_{2r-1}<B^-_{2r-1}\},\\
X^r_2&=&\{2\lambda: 2\lambda_{2r-1}\geq B^-_{2r-1}\},
\end{eqnarray*}
and for  $1\leq s\leq \frac{n_0}{2}+1$,
 \begin{eqnarray*}
Y^s_1&=&\{2\lambda: 2\lambda_{2s}\leq B^+_{2s}\},\\
Y^s_2&=&\{2\lambda: 2\lambda_{2s}> B^+_{2s}\}.
\end{eqnarray*}
Now we let
\begin{equation}\label{sn}
S(n)=\{(i,j):\ 0\leq i\leq \frac{n_1+1}{2}, \ 0\leq j\leq
\frac{n_0}{2}\}.
\end{equation} Then for $(i,j)\in S(n)$, we form the
intersection
\begin{eqnarray*}
L_{i,j}(n,\sigma,\alpha)= &&(X^1_2\cap\cdots\cap X^{i}_2\cap X^{i+1}_1\cap\cdots\cap X^{\frac{n_1+1}{2}+1}_1)\cap \\
&&(Y^1_2\cap\cdots\cap Y^{j}_2\cap Y^{j+1}_1\cap\cdots\cap Y^{\frac{n_0}{2}+1}_1).
\end{eqnarray*}
More precisely, $2\lambda\in L_{i,j}(n,\sigma,\alpha)$ if and only if
\begin{equation}
2\lambda_{2i-1}\geq  B^-_{2i-1}\geq 2\lambda_{2i+1}\ \ \mbox{and}\ \
2\lambda_{2j}\geq B^+_{2j+2}\geq 2\lambda_{2j+2},
\end{equation}
where
\[B^-_{2i-1}=\sigma-\frac{n+\alpha+1}{2}+2i-1+1=\sigma-\frac{n+\alpha}{2}+2i-\frac{1}{2}\] and
\[B^+_{2j+2}=-\sigma-\frac{n+\alpha+1}{2}+2j+2-1=-\sigma-\frac{n+\alpha}{2}+2j+\frac{1}{2}.\]
Note that $L_{a,b}(n,\sigma,\alpha)$ may be empty. If it is nonempty, then we shall identify it with the direct sum of all the $K$-representations
$V_{2\lambda+\frac{\alpha}{2}{\mbox{\rm\bf 1}}}$ with $2\lambda\in L_{a,b}(n,\sigma,\alpha)$.

\vspace{0.2in}\noindent{\bf Remark.}  Let $k=[n/2]$. Then
\begin{equation}
S(n)=\left\{\begin{array}{ll} \{(i,j):\ 0\leq
i\leq k+1,\ 0\leq j\leq k\}&\mbox{$n$ odd,}\\
&\\
\{(i,j):\ 0\leq i,j\leq k\}&\mbox{$n$ even.}
\end{array}\right.
\end{equation}

\vspace{0.2in}
The following combinatorial result is elementary.
\vspace{0.2in}\begin{lemma}
\begin{enumerate}
\item[(i)] Assume that $(i,j)\in S(n)$ and  $i-j\leq -1$. Then $L_{i,j}\neq\emptyset$ if and only if $i-j\geq-\sigma+\frac{1}{2}$. In particular, if $\sigma\leq -1/2$, then all such $L_{i,j}$ is empty.
\item[(ii)] For $0\leq i\leq k$, $L_{i,i}$ is always nonempty.
\item[(iii)] If $(j+1,j)\in S(n)$, then $L_{j+1,j}$ is always nonempty.
\item[(iv)]  Assume that $(i,j)\in S(n)$ and $i-j\geq 2$. Then $L_{i,j}\neq\emptyset$ if and only if $i-j\leq-\sigma+\frac{1}{2}$. In particular, if $\sigma\geq 1/2$, then all such $L_{i,j}$ is empty.
\end{enumerate}
\end{lemma}

\vspace{0.2in}
\begin{thm} \label{evenodd}
{\rm (Case 2a)}
Assume that $n+\alpha$ is even, $\sigma+{\textstyle \frac{1}{2}}\in{\mathbb Z}$, and $\tilde{\sigma} =\sigma+ \frac{n+\alpha+1}{2}$
is odd.

\vspace{0.05in}
\begin{enumerate}
\item[(i)] If $\sigma\geq 1/2$, then
\[{I^\alpha(\sigma)}=\bigoplus\{L_{i,j}(n,\sigma,\alpha):\ (i,j)\in S(n), \ -1\leq j-i\leq r_1\},\]
where
\[r_1=\min\left(\sigma-{\textstyle \frac{1}{2}},{\textstyle \left[\frac{n}{2}\right]}\right),\]
and each $L_{i,j}(n,\sigma,\alpha)$ which appears in the sum forms an
irreducible constituent of ${I^\alpha(\sigma)}$. The module diagram of ${I^\alpha(\sigma)}$ can
be obtained from Figure 3 by removing those $L_{a,b}(n,\sigma,\alpha)$
which are empty. In particular, the socle series of ${I^\alpha(\sigma)}$ is given
by
\[\mbox{\rm Soc}^l({I^\alpha(\sigma)}=\left\{
\begin{array}{ll}
\bigoplus_{(i,j)\in S(n),\ -1\leq j-i\leq l-2}L_{i,j}(n,\sigma,\alpha)&1\leq l\leq r_1+1,\\
{I^\alpha(\sigma)}&l\geq r_1+2.
\end{array}\right.\]
An irreducible constituent $L_{i,j}(n,\sigma,\alpha)$ of ${I^\alpha(\sigma)}$ is
unitarizable if and only if $i=j+1$ or $\frac{1}{2}\leq \sigma\leq
[\frac{n}{2}]+{\textstyle \frac{1}{2}}$ and $j-i=r_1$.

\vspace{0.2in}
\item[(ii)] If $\sigma\leq -1/2$, then
\[{I^\alpha(\sigma)}=\bigoplus\{L_{i,j}(n,\sigma,\alpha):\ (i,j)\in S(n),\ 0\leq i-j\leq r_2\},\]
where
\[r_2=\min(-\sigma+{\textstyle \frac{1}{2}}, {\textstyle \left[\frac{n+1}{2}\right]}),\]
and each $L_{i,j}(n,\sigma,\alpha)$ which appears in the sum forms an
irreducible constituent of ${I^\alpha(\sigma)}$. The module diagram of ${I^\alpha(\sigma)}$ can
be obtained from Figure 3 by removing those $L_{a,b}(n,\sigma,\alpha)$
which are empty. In particular, the socle series of ${I^\alpha(\sigma)}$ is given
by
\[\mbox{\rm Soc}^l({I^\alpha(\sigma)})=\left\{\begin{array}{ll}
\bigoplus _{(i,j)\in S(n),\ r_2-l+1\leq i-j\leq r_2}L_{i,j}(n,\sigma,\alpha) &1\leq l\leq r_2,\\
{I^\alpha(\sigma)}&l\geq r_2+1.
\end{array}
\right.\] An irreducible constituent $L_{i,j}(n,\sigma,\alpha)$ of
${I^\alpha(\sigma)}$ is unitarizable if and only if $i=j$ or $-
\left[\frac{n+1}{2}\right]+{\textstyle \frac{1}{2}}\leq\sigma\leq -{\textstyle \frac{1}{2}}$ and $i-j=r_2$.
\end{enumerate}
\end{thm}

\vspace{0.2in} The module diagram of ${I^\alpha(\sigma)}$ in the case when $n+\alpha$ is
even and $\sigma+  \frac{n+\alpha+1}{2}$ is odd can be obtained from
the following rectangle (Figure 3) by removing those $L_{a,b}(n,\sigma,\alpha)$
which are empty. Here $k=[n/2]$ and we write a constituent
$L_{a,b}(n,\sigma,\alpha)$ simply as $L_{a,b}$. Note also that when $n$
is even, the spaces $L_{k+1,j}$ ($0\leq j\leq k$) are not defined. In
this case, the rectangle below reduces to a square.

{\tiny
 \vspace{0.2in}
\begin{center}
\setlength{\unitlength}{0.8mm}

\begin{picture}(160,140)(-20,-10)

\put(60,125){\makebox(10,10){$  L_{0,k} $}}
\put(45,110){\makebox(10,10){$  L_{1,k} $}}
\put(75,110){\makebox(10,10){$  L_{0,k-1} $}}
\put(30,95){\makebox(10,10){$L_{2,k} $}}
\put(60,95){\makebox(10,10){$L_{1,k-1} $}}
 \put(-15,50){\makebox(10,10){$L_{k+1,k} $}}
\put(-7,60){\line(1,1){5}}
 \put(15,50){\makebox(10,10){$L_{k,k-1} $}}
\put(75,50){\makebox(10,10){$L_{2,1} $}}
\put(105,50){\makebox(10,10){$L_{1,0} $}}
\put(30,5){\makebox(10,10){$L_{k+1,1} $}} \put(16,29){.}
\put(19,26){.} \put(22,23){.} \put(32,13){\line(-1,1){5}}
\put(45,-10){\makebox(10,10){$L_{k+1,0} $}}
\put(13,32){\line(-1,1){5}} \put(0,65){\makebox(10,10){$L_{k,k}
$}} \put(30,65){\makebox(10,10){$L_{k-1,k-1} $}}
\put(90,65){\makebox(10,10){$L_{1,1} $}}
\put(105,80){\makebox(10,10){$L_{0,1} $}}
\put(120,65){\makebox(10,10){$L_{0,0} $}}
 \put(0,35){\makebox(10,10){$L_{k+1, k-1} $}}
\put(90,35){\makebox(10,10){$L_{2,0} $}}
\put(45,20){\makebox(10,10){$L_{k,1} $}}
\put(60,5){\makebox(10,10){$L_{k,0} $}}
\put(38,7){\line(1,-1){10}} \put(62,7){\line(-1,-1){10}}
\put(53,118){\line(1,1){9}} \put(77,118){\line(-1,1){9}}
\put(38,103){\line(1,1){9}} \put(62,103){\line(-1,1){9}}
\put(68,103){\line(1,1){9}} \put(88,107){\line(-1,1){5}}
 \put(89,104){.}
\put(92,101){.} \put(95,98){.} \put(32,97){\line(-1,-1){5}}
\put(38,97){\line(1,-1){5}} \put(62,97){\line(-1,-1){5}}
\put(68,97){\line(1,-1){5}} \put(92,97){\line(-1,-1){5}}
\put(98,97){\line(1,-1){5}} \put(8,75){\line(1,1){5}}
\put(32,75){\line(-1,1){5}} \put(38,75){\line(1,1){5}}
\put(92,75){\line(-1,1){5}} \put(98,75){\line(1,1){5}}
\put(122,75){\line(-1,1){5}} \put(8,67){\line(1,-1){5}}
\put(32,67){\line(-1,-1){5}} \put(38,67){\line(1,-1){5}}
\put(92,67){\line(-1,-1){5}} \put(98,67){\line(1,-1){5}}
\put(122,67){\line(-1,-1){5}} \put(2,43){\line(-1,1){5}}
\put(8,43){\line(1,1){5}} \put(32,43){\line(-1,1){5}}
\put(38,43){\line(1,1){5}} \put(68,43){\line(1,1){5}}
\put(92,43){\line(-1,1){5}} \put(98,43){\line(1,1){5}}
\put(32,37){\line(-1,-1){9}} \put(38,37){\line(1,-1){9}}
\put(58,33){\line(-1,-1){5}} \put(38,37){\line(1,-1){9}}
\put(72,33){\line(1,-1){5}} \put(92,37){\line(-1,-1){5}}
 \put(47,22){\line(-1,-1){10}}
\put(53,22){\line(1,-1){10}} \put(72,17){\line(-1,-1){5}}
\put(77,22){.} \put(80,25){.}
 \put(83,28){.}
\put(17,82){.} \put(20,85){.} \put(23,88){.}
 \put(55,70){.}
\put(59,70){.} \put(63,70){.} \put(67,70){.}

 \end{picture}

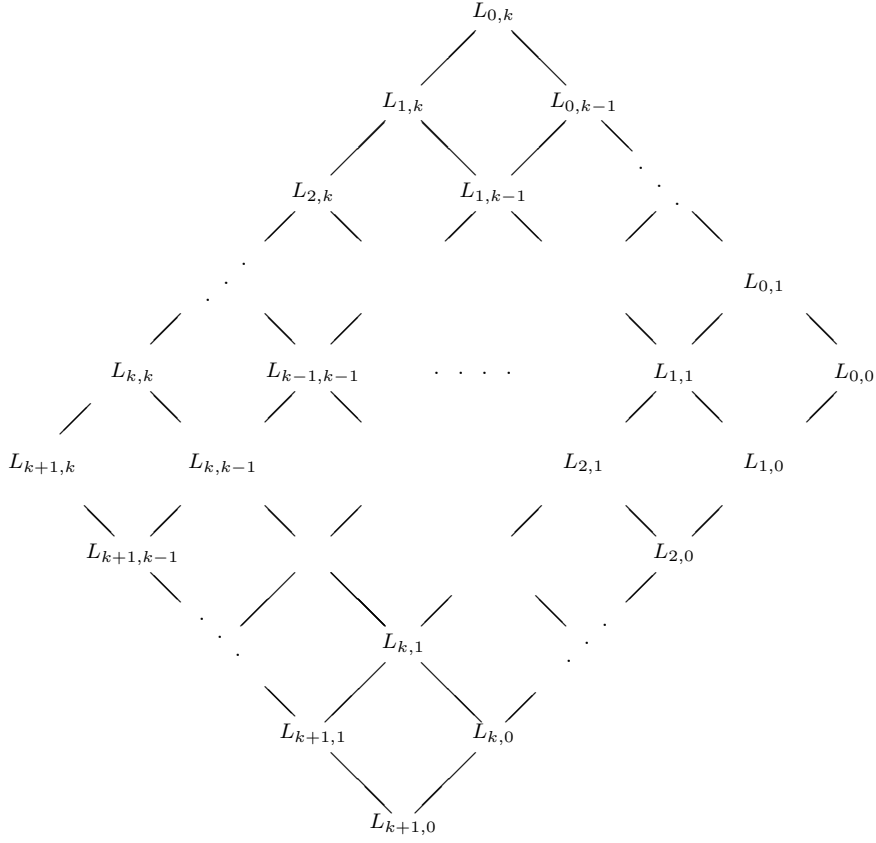
Figure 3: Diagram for ${I^\alpha(\sigma)}$ when $n+\alpha$ is even and $\sigma+
\frac{n+\alpha+1}{2}$ is odd

\end{center}}

\vspace{0.2in}
\vspace{0.2in}\noindent{\bf Case 2b:}  $\sigma+\frac{n+\alpha+1}{2}$ is even

As in Case 2a, exactly one barrier is effective along each coordinate axis. But now $\ell^+_j$ is
effective when $j$ is odd, and $\ell^-_j$ is effective when $j$ is
even. So we have to define $L_{a,b}(n,\sigma,\alpha)$ differently.
For $1\leq r\leq \frac{n_1+1}{2}+1$, let
\begin{eqnarray*}
X^r_1&=&\{2\lambda: 2\lambda_{2r-1}\leq B^+_{2r-1}\},\\
X^r_2&=&\{2\lambda: 2\lambda_{2r-1}> B^-_{2r-1}\},
\end{eqnarray*}
and for  $1\leq s\leq \frac{n_0}{2}+1$,
 \begin{eqnarray*}
Y^s_1&=&\{2\lambda: 2\lambda_{2s}< B^-_{2s}\},\\
Y^s_2&=&\{2\lambda: 2\lambda_{2s}\geq B^-_{2s}\}.
\end{eqnarray*}
We define the set $S(n)$, as in Case 2a in equation \eqref{sn}. Then for $(i,j)\in S(n)$, we form the
intersection
\begin{eqnarray*}
L_{i,j}(n,\sigma,\alpha)= & & (X^1_2\cap\cdots\cap X^{i}_2\cap X^{i+1}_1\cap\cdots\cap X^{\frac{n_1+1}{2}+1}_1)\cap \\
& & (Y^1_2\cap\cdots\cap Y^{j}_2\cap Y^{j+1}_1\cap\cdots\cap Y^{\frac{n_0}{2}+1}_1).\end{eqnarray*}
More precisely, $2\lambda\in L_{i,j}(\sigma,\alpha)$ if and only if
\begin{equation}
2\lambda_{2i-1}\geq  B^+_{2i+1}\geq 2\lambda_{2i+1}\ \ \mbox{and}\ \
2\lambda_{2j}\geq B^-_{2j}\geq 2\lambda_{2j+2}.
 \end{equation}
where
\[B^+_{2i+1}=-\sigma-\frac{n+\alpha+1}{2}+2i+1-1=-\sigma-\frac{n+\alpha}{2}+2i-\frac{1}{2}\] and
\[B^-_{2j}=\sigma-\frac{n+\alpha+1}{2}+2j+1=\sigma-\frac{n+\alpha}{2}+2j+\frac{1}{2}.\]

\vspace{0.2in}
\begin{thm} \label{eveneven}
{\rm (Case 2b)}
Assume that $n+\alpha$ is even, $\sigma+{\textstyle \frac{1}{2}}\in{\mathbb Z}$, and $\tilde{\sigma} =\sigma+ \frac{n+\alpha+1}{2}$
is even.

\vspace{0.05in}
\begin{enumerate}
\item[(i)] If $\sigma\geq 1/2$, then
\[{I^\alpha(\sigma)}=\bigoplus\{L_{i,j}(n,\sigma,\alpha):\ (i,j)\in S(n), \ 0\leq i-j\leq r_2\},\]
where
\[r_2=\min(\sigma+{\textstyle \frac{1}{2}}, {\textstyle \left[\frac{n+1}{2}\right]}),\]
and each $L_{i,j}(n,\sigma,\alpha)$ which appears in the sum forms an
irreducible constituent of ${I^\alpha(\sigma)}$. The module diagram of ${I^\alpha(\sigma)}$ can
be obtained from Figure 4 by removing those $L_{a,b}(n,\sigma,\alpha)$
which are empty. In particular, the socle series of ${I^\alpha(\sigma)}$ is given
by
\[\mbox{\rm Soc}^l({I^\alpha(\sigma)}=\left\{
\begin{array}{ll}
\bigoplus _{(i,j)\in S(n),
 0\leq i-j\leq l-1} L_{i,j}(n,\sigma,\alpha)
 &1\leq l\leq r_2,\\
 &\\
{I^\alpha(\sigma)}&l\geq r_2+1.
\end{array}\right.\]
 An irreducible constituent $L_{i,j}(n,\sigma,\alpha)$ of ${I^\alpha(\sigma)}$ is unitarizable
if and only if $i=j$ or ${\textstyle \frac{1}{2}}\leq\sigma\leq \left[\frac{n+1}{2}\right]-{\textstyle \frac{1}{2}}$ and $i-j=r_2$.

\vspace{0.2in}
\item[(ii)] If $\sigma\leq -1/2$, then
\[{I^\alpha(\sigma)}=\bigoplus\{L_{i,j}(n,\sigma,\alpha):\ (i,j)\in S(n),\ -1\leq j-i\leq r_1\},\]
where
\[r_1=\min\left(-\sigma-{\textstyle \frac{1}{2}},{\textstyle \left[\frac{n}{2}\right]}\right),\]
and each $L_{i,j}(n,\sigma,\alpha)$ which appears in the sum forms an
irreducible constituent of ${I^\alpha(\sigma)}$. The module diagram of ${I^\alpha(\sigma)}$ can
be obtained from Figure 4 by removing those $L_{a,b}(n,\sigma,\alpha)$
which are empty. In particular, the socle series of ${I^\alpha(\sigma)}$ is given
by
\[\mbox{\rm Soc}^l({I^\alpha(\sigma)})=\left\{\begin{array}{ll}
\bigoplus _{(i,j)\in S(n),\ r_1-l+1\leq j-i\leq r_1} L_{i,j}(n,\sigma,\alpha)
&1\leq l\leq r_1+1,\\
&\\
{I^\alpha(\sigma)}&l\geq r_1+2.
\end{array}
\right.\]
An irreducible constituent $L_{i,j}(n,\sigma,\alpha)$ of ${I^\alpha(\sigma)}$
is unitarizable if and only if $i=j+1$ or $-[\frac{n}{2}]-{\textstyle \frac{1}{2}}\leq
\sigma\leq -{\textstyle \frac{1}{2}}$ and $j-i=r_1$.

\end{enumerate}
\end{thm}

\vspace{0.2in}  Finally, we describe the module diagram  of ${I^\alpha(\sigma)}$ in the case
when both $n+\alpha$   and $\sigma+  \frac{n+\alpha+1}{2}$ are even.
  It can be obtained from the following rectangle (Figure 4) by removing
those $L_{a,b}(n,\sigma,\alpha)$ which are empty.   As before,
$k=[n/2]$, and  we write a constituent $L_{a,b}(n,\sigma,\alpha)$
simply as $L_{a,b}$. When $n$ is even, $L_{k+1, j}$ ($0\leq j\leq k$) is empty so that this rectangle reduces to a
square.

\vspace{0.2in}{\tiny
\begin{center}
\setlength{\unitlength}{0.8mm}

\begin{picture}(160,140)(-20,-10)

\put(60,125){\makebox(10,10){$L_{k+1,0} $}}
\put(45,110){\makebox(10,10){$L_{k,0} $}}
\put(75,110){\makebox(10,10){$L_{k+1,1} $}}
\put(30,95){\makebox(10,10){$L_{k-1,0} $}}
\put(60,95){\makebox(10,10){$L_{k,1} $}}

 \put(-15,50){\makebox(10,10){$L_{0,0}$}}
\put(-7,60){\line(1,1){5}}
 \put(15,50){\makebox(10,10){$L_{1,1} $}}
\put(75,50){\makebox(10,10){$L_{k-1,k-1} $}}
\put(105,50){\makebox(10,10){$L_{k,k} $}}

\put(30,5){\makebox(10,10){$L_{0,k-1} $}}

\put(90,35){\makebox(10,10){$L_{k-1,k} $}}

\put(16,29){.} \put(19,26){.} \put(22,23){.}
\put(32,13){\line(-1,1){5}} \put(45,-10){\makebox(10,10){$L_{0,k}
$}} \put(13,32){\line(-1,1){5}}

\put(0,65){\makebox(10,10){$L_{1,0} $}}
\put(30,65){\makebox(10,10){$L_{2,1} $}}
\put(90,65){\makebox(10,10){$L_{k,k-1} $}}
\put(120,65){\makebox(10,10){$L_{k+1,k} $}}

\put(105,80){\makebox(10,10){$L_{k+1,k-1} $}}

 \put(0,35){\makebox(10,10){$L_{0,1} $}}

\put(45,20){\makebox(10,10){$L_{1,k-1} $}}

\put(60,5){\makebox(10,10){$L_{1,k} $}}
\put(38,7){\line(1,-1){10}} \put(62,7){\line(-1,-1){10}}
\put(53,118){\line(1,1){9}} \put(77,118){\line(-1,1){9}}
\put(38,103){\line(1,1){9}} \put(62,103){\line(-1,1){9}}
\put(68,103){\line(1,1){9}} \put(88,107){\line(-1,1){5}}
 \put(89,104){.}
\put(92,101){.} \put(95,98){.} \put(32,97){\line(-1,-1){5}}
\put(38,97){\line(1,-1){5}} \put(62,97){\line(-1,-1){5}}
\put(68,97){\line(1,-1){5}} \put(92,97){\line(-1,-1){5}}
\put(98,97){\line(1,-1){5}} \put(8,75){\line(1,1){5}}
\put(32,75){\line(-1,1){5}} \put(38,75){\line(1,1){5}}
\put(92,75){\line(-1,1){5}} \put(98,75){\line(1,1){5}}
\put(122,75){\line(-1,1){5}} \put(8,67){\line(1,-1){5}}
\put(32,67){\line(-1,-1){5}} \put(38,67){\line(1,-1){5}}
\put(92,67){\line(-1,-1){5}} \put(98,67){\line(1,-1){5}}
\put(122,67){\line(-1,-1){5}} \put(2,43){\line(-1,1){5}}
\put(8,43){\line(1,1){5}} \put(32,43){\line(-1,1){5}}
\put(38,43){\line(1,1){5}} \put(68,43){\line(1,1){5}}
\put(92,43){\line(-1,1){5}} \put(98,43){\line(1,1){5}}
\put(32,37){\line(-1,-1){9}} \put(38,37){\line(1,-1){9}}
\put(58,33){\line(-1,-1){5}} \put(38,37){\line(1,-1){9}}
\put(72,33){\line(1,-1){5}} \put(92,37){\line(-1,-1){5}}
 \put(47,22){\line(-1,-1){10}}
\put(53,22){\line(1,-1){10}} \put(72,17){\line(-1,-1){5}}
\put(77,22){.} \put(80,25){.}
 \put(83,28){.}
\put(17,82){.} \put(20,85){.} \put(23,88){.}
 \put(55,70){.}
\put(59,70){.} \put(63,70){.} \put(67,70){.}
 \end{picture}

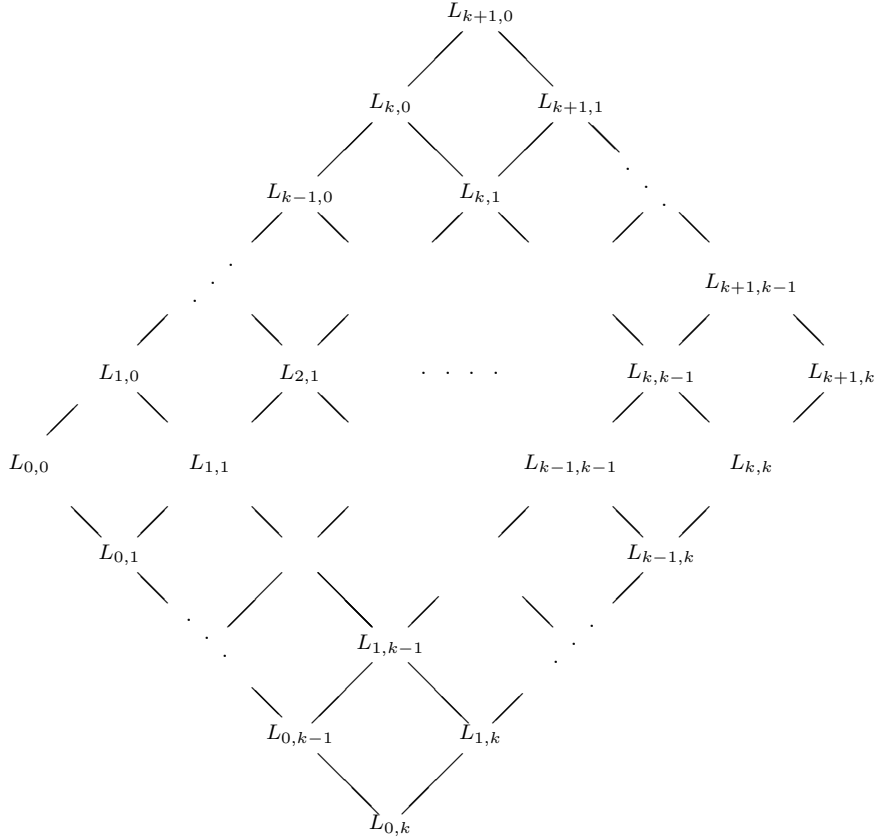
Figure 4: Diagram for ${I^\alpha(\sigma)}$ when both $n+\alpha$   and $\sigma+
\frac{n+\alpha+1}{2}$ are even

\end{center}}

\section{The image of $\psi_{p,q}$ in ${I^\alpha(\sigma)}$}
\label{Image}

In this section, we will describe the image of $\psi _{p,q}$ in ${I^\alpha(\sigma)}$:
\[\psi _{p,q}: \ \ \Omega ^{p,q}\hookrightarrow I^\alpha(\sigma),\]
where
\begin{equation*}
\sigma = \frac{p+q}{2}-\frac{n+1}{2}, \ \ \text{and} \ \ \alpha \equiv p-q \ (\mathrm{mod}\  4).
\end{equation*}
We allow $(p,q)=(0,0)$, in which case we understand $\Omega ^{0,0}$ to be the trivial representation of $G={\widetilde{\mathrm{Sp}}(2n,{\mathbb R})}$.

Denote $m=p+q$. Observe that
\[n+\alpha \equiv n+m \ (\mathrm{mod}\  2),\]
and
\[\tilde{\sigma} = \sigma+\frac{n+\alpha+1}{2}=\frac{m+\alpha}{2}\equiv p\ (\mathrm{mod}\  2).\]

Recall that our analysis of ${I^\alpha(\sigma)}$ is divided into 4 cases:

\vspace{0.2in}
\begin{center}
\begin{tabular}{|c|c|c|}\hline
  $\tilde{\sigma} \backslash n+\alpha $& odd & even \\ \hline odd & Case 1a & Case 2a\\
  \hline even &
Case 1b   & Case 2b\\  \hline
\end{tabular}
\end{center}

Correspondingly, our analysis of $\Omega ^{p,q}$ will also be divided
into 4 cases:

\vspace{0.2in}
\begin{center}
\begin{tabular}{|c|c|c|}\hline
  $p \backslash n+m$& odd & even \\ \hline odd & Case 1a & Case 2a\\
  \hline even &
Case 1b   & Case 2b\\  \hline
\end{tabular}
\end{center}

\vspace{0.2in}

Throughout this section, we assume that we are at a point of reducibility, namely
\[\sigma +\frac{n+\alpha +1}{2}\in {\mathbb Z}.\]
If $\frac{p+q-(n+1)}{2}=\sigma$ and $p-q\equiv \alpha$ (mod $4$), we say that we have a possible embedding of $\Omega ^{p,q}$ (into
$I^\alpha(\sigma)$).

We shall give the detailed description of $\Omega ^{p,q}$ in ${I^\alpha(\sigma)}$ in the following two subsections. We end this section with a corollary which gives conceptual underpinning to our results.

\begin{thm} \label{summary} The relationship between $\Omega ^{p,q}$'s and ${I^\alpha(\sigma)}$ is as follows.
\begin{itemize}
\item[(a)] $-\rho_n\leq \sigma <0$: The irreducible submodules of $I^\alpha(\sigma)$ are given by the possible embeddings of $\Omega ^{p,q}$'s, and all of them are unitary.
\item[(b)] Unitary axis ($\sigma=0$) (when $n+\alpha$ is odd):
\[I^\alpha(0)=\bigoplus_{\begin{smallmatrix}p+q=n+1\\
p-q\equiv \alpha \ (\mathrm{mod}\  4)\end{smallmatrix}}\Omega ^{p,q}.\]
\item[(c)] $\sigma >0$: The reducibilities of $I^\alpha(\sigma)$ are completely accounted for by the possible embeddings of $\Omega ^{p,q}$'s.
\begin{itemize}
\item See Part b) of Theorems \ref{Ioo}, \ref{Ioe}, \ref{Ieo} and \ref{Iee} for the precise statements.
\end{itemize}
\end{itemize}
\end{thm}

\subsection{The image of $\psi_{p,q}$ in ${I^\alpha(\sigma)}$: $m$ and $n$ different parity}

\[
\]

\noindent{\bf Case 1: $n+m$ odd.}

\vspace{0.2in} \noindent We consider two subcases:
\[\mbox{{\em Case 1a:}  $p$ is odd\hspace{0.5in} {\em Case 1b:} $p$ is even.}\]

\noindent {\bf Case 1a}: we write
\[p=2i+1, \ \ q=2j+\epsilon,\]
so that $p+q=m$. Here
\[\epsilon=\left\{\begin{array}{ll}
1 &\mbox{$m$ even,}\\
0 &\mbox{$m$ odd.}
\end{array}\right.\]
Note that
\[i+j=\frac{m-1-\epsilon}{2}=\frac{m+n_0-(n+1)}{2}=\frac{n_0}{2}+\sigma.\]

\noindent {\bf Case 1b}: we write
\[p=2i, \ \ q=2j+\epsilon ',\]
so that $p+q=m$. Here
\[\epsilon '=\left\{\begin{array}{ll}
0 &\mbox{$m$ even,}\\
1 &\mbox{$m$ odd.}
\end{array}\right.\]
Note that
\[i+j=\frac{m-\epsilon '}{2}=\frac{m+(n_1+1)-(n+1)}{2}=\frac{n_1+1}{2}+\sigma.\]

In accordance with \ref{defk}, we let
\begin{equation*}
%\label{defk2}
k=\left\{\begin{array}{ll}\displaystyle
\frac{n_0}{2}=[\frac{n}{2}]&\mbox{\hspace{0.5in} if $p$ is odd,}\\
&\\  \displaystyle \frac{n_1+1}{2}=[\frac{n+1}{2}]&\mbox{\hspace{0.5in} if $p$ is even.}
\end{array}\right.
\end{equation*}

Recall that we have the following module diagrams:
\vspace{0.2in} {\tiny \setlength{\unitlength}{0.8mm}

\begin{center}

\begin{picture}(120,40)(0,30)
%\begin{picture}(120,65)(0,5)

\put(0,60){\makebox(10,10){$\scriptstyle R_{k,0} $}}
\put(20,60){\makebox(10,10){$\scriptstyle R_{k-1,1}   $}}

\put(40,65){$\cdot$} \put(50,65){$\cdot$} \put(60,65){$\cdot$}
\put(70,65){$\cdot$} \put(80,65){$\cdot$}

\put(90,60){\makebox(10,10){$\scriptstyle R_{1,k-1}  $}}
\put(110,60){\makebox(10,10){$\scriptstyle R_{0,k}  $}}

\put(25,35){\makebox(10,10){$\scriptstyle R_{r_1,0} $}}
\put(55,35){\makebox(10,10){$\scriptstyle ... $}}
\put(85,35){\makebox(10,10){$\scriptstyle R_{0,r_1} $}}

\put(27,44){\line(-1,1){5}} \put(33,44){\line(1,1){5}}
\put(57,44){\line(-1,1){5}} \put(63,44){\line(1,1){5}}
\put(87,44){\line(-1,1){5}} \put(93,44){\line(1,1){5}}

\put(20,51){.} \put(18,53){.} \put(16,55){.}

\put(14,57){\line(-1,1){5}} \put(16,57){\line(1,1){5}}
\put(30,57){\line(-1,1){5}}

\put(100,51){.} \put(102,53){.} \put(104,55){.}
\put(106,57){\line(1,1){5}} \put(104,57){\line(-1,1){5}}
\put(90,57){\line(1,1){5}}

\end{picture}

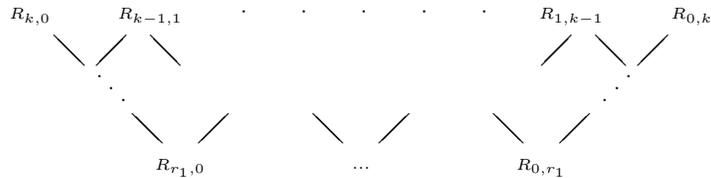
Case 1: Module diagram for ${I^\alpha(\sigma)}$ $(-k\leq \sigma\leq 0)$; $r_1=k+\sigma$
\end{center}}

\vspace{0.2in} {\tiny \setlength{\unitlength}{0.8mm}

\begin{center}

\begin{picture}(120,40)(0,5)
%\begin{picture}(120,65)(0,5)

\put(0,10){\makebox(10,10){$\scriptstyle R_{k,0}   $}}
\put(20,10){\makebox(10,10){$\scriptstyle R_{k-1,1}$}}

\put(40,15){$\cdot$} \put(50,15){$\cdot$} \put(60,15){$\cdot$}
\put(70,15){$\cdot$} \put(80,15){$\cdot$}

\put(90,10){\makebox(10,10){$\scriptstyle R_{1,k-1}$}}
\put(110,10){\makebox(10,10){$\scriptstyle R_{0,k}$}}

\put(25,35){\makebox(10,10){$\scriptstyle R_{r_2,0}$}}
\put(53,35){\makebox(10,10){$...$}}
\put(85,35){\makebox(10,10){$\scriptstyle R_{0,r_2}$}}

\put(24,37){\line(-1,-1){5}} \put(29,37){\line(1,-1){5}}
\put(55,37){\line(-1,-1){5}} \put(60,37){\line(1,-1){5}}
\put(85,37){\line(-1,-1){5}} \put(91,37){\line(1,-1){5}}

\put(13,26){.} \put(15,28){.} \put(17,30){.}
\put(102,26){.} \put(100,28){.} \put(98,30){.}

\put(6,19){\line(1,1){5}} \put(109,19){\line(-1,1){5}}
\put(23,19){\line(-1,1){5}} \put(25,19){\line(1,1){5}}
\put(92,19){\line(-1,1){5}} \put(94,19){\line(1,1){5}}

\end{picture}

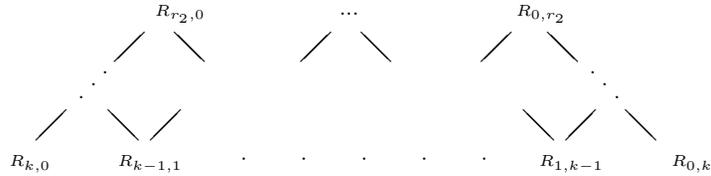
Case 1: Module diagram for ${I^\alpha(\sigma)}$ ($1\leq \sigma \leq k$); $r_2=k-\sigma$

\end{center}}

\vspace{0.2in}

We introduce one notation. For an irreducible subquotient $R_{s,t}$ of ${I^\alpha(\sigma)}$ where $\sigma >0$, denote by $\prec R_{s,t}\succ$ the submodule of ${I^\alpha(\sigma)}$ generated by $R_{s,t}$:
\begin{equation}
\prec R_{s,t}\succ =\oplus \{R_{i,j}: i\geq s, j\geq t\}.
\end{equation}

Recall that
\[\sigma =\frac{m-(n+1)}{2}, \ \ \text{and} \ \ p-q\equiv \alpha \ (\text{mod}\ 4). \]

\vspace{0.2in}\begin{thm} \label{Ioo}
{\rm (Case 1a)}
Assume that $n+m$ is odd, and $p$ is odd.
\begin{enumerate}
\item[(a)] If $\sigma\leq 0$, then
\[\Omega ^{p,q}=R_{\frac{p-1}{2},\frac{q-\epsilon}{2}}.\]
\vspace{0.2in}
\item[(b)] If $\sigma\geq 1$,
then
\[\Omega ^{p,q}=\prec R_{s,t}\succ,\]
where $s =\max\left(0, \frac{n-q}{2}\right)$, and $t =\max\left(0, \frac{n+1-\epsilon-p}{2}\right)$.
\end{enumerate}
\end{thm}

\vspace{0.2in}\noindent{\bf Remarks:} \begin{enumerate} \item[(a)] For $(i,j)=(\frac{p-1}{2},\frac{q-\epsilon }{2})$, we have noted that $i+j=\frac{n_0}{2}+\sigma =r_1$. Thus $R_{i,j}$ is at the bottom layer of the module diagram of ${I^\alpha(\sigma)}$, namely is an irreducible submodule. The collection of $R_{i,j}$'s exhausts the set of all irreducible submodules of ${I^\alpha(\sigma)}$.
\item[(b)] If $q\leq n$, and $p\leq n+1-\epsilon $, then $s=\max\left(0, \frac{n-q}{2}\right)=\frac{n-q}{2}$, $t =\max\left(0, \frac{n+1-\epsilon -p}{2}\right)= \frac{n+1-\epsilon -p}{2}$, and $s+t=\frac{n-\epsilon }{2}-\sigma =\frac{n_0}{2}-\sigma = r_2$. Such $R_{s,t}$'s are exactly those at the top layer of the module diagram of ${I^\alpha(\sigma)}$, namely irreducible quotients. The rest of $R_{s,t}$'s are those on the ``left boundary" ($s=0$) or the ``right boundary" ($t=0$).
\end{enumerate}

\vspace{0.2in}\begin{thm} \label{Ioe}
{\rm (Case 1b)}
Assume that $n+m$ is odd, and $p$ is even.
\begin{enumerate}
\item[(a)] If $\sigma\leq 0$, then
\[\Omega ^{p,q}=R_{\frac{p}{2},\frac{q-\epsilon '}{2}}.\]
\vspace{0.2in}
\item[(b)] If $\sigma\geq 1$,
then
\[\Omega ^{p,q}=\prec R_{s,t}\succ,\]
where $s =\max\left(0, \frac{n+1-q}{2}\right)$, and $t =\max\left(0, \frac{n+1-\epsilon '-p}{2}\right)$.
\end{enumerate}
\end{thm}

\vspace{0.2in}\noindent{\bf Remarks:} \begin{enumerate} \item[(a)] For $(i,j)=(\frac{p}{2},\frac{q-\epsilon '}{2})$, we have noted that $i+j=\frac{n_1+1}{2}+\sigma =r_1$. Thus $R_{i,j}$ is at the bottom layer of the module diagram of ${I^\alpha(\sigma)}$, namely is an irreducible submodule. The collection of $R_{i,j}$'s exhausts the set of all irreducible submodules of ${I^\alpha(\sigma)}$.
\item[(b)] If $q\leq n+1$, and $p\leq n+1-\epsilon '$, then $s=\max\left(0, \frac{n+1-q}{2}\right)=\frac{n+1-q}{2}$, $t =\max\left(0, \frac{n+1-\epsilon '-p}{2}\right)= \frac{n+1-\epsilon '-p}{2}$, and $s+t=\frac{n+1-\epsilon '}{2}-\sigma =\frac{n_1+1}{2}-\sigma = r_2$. Such $R_{s,t}$'s are exactly those at the top layer of the module diagram of ${I^\alpha(\sigma)}$, namely irreducible quotients. The rest of $R_{s,t}$'s are those on the ``left boundary" ($s=0$) or the ``right boundary" ($t=0$).
\end{enumerate}

\vspace{0.2in}

\subsection{The image of $\psi_{p,q}$ in ${I^\alpha(\sigma)}$: $m$ and $n$ same parity}
\[
\]
\noindent{\bf Case 2: $n+m$ even.}

\vspace{0.2in} \noindent We consider two subcases:
\[\mbox{{\em Case 2a:}  $p$ is odd\hspace{0.5in} {\em Case 2b:} $p$ is even.}\]

We introduce one similar notation in this case. For an irreducible subquotient $L_{s,t}$ of ${I^\alpha(\sigma)}$ where $\sigma >0$, denote by $\prec L_{s,t}\succ$ the submodule of ${I^\alpha(\sigma)}$ generated by $L_{s,t}$:
\begin{equation}
\prec L_{s,t}\succ =\begin{cases} \oplus \{L_{i,j}: i\geq s, j\leq t\}, \ \ \text{(Case 2a)},\\
\oplus \{L_{i,j}: i\leq s, j\geq t\}, \ \ \text{(Case 2b)}.
\end{cases}
\end{equation}

As in Subsection \ref{Subeven}, we let $k=[n/2]$. First recall the module diagram of ${I^\alpha(\sigma)}$ for Case 2a:

{\tiny
 \vspace{0.2in}
\begin{center}
\setlength{\unitlength}{0.8mm}

\begin{picture}(160,60)(-20,10)
%\begin{picture}(160,140)(-20,-10)

\put(-15,50){\makebox(10,10){$L_{k+1,k} $}}
\put(-7,60){\line(1,1){5}}
\put(15,50){\makebox(10,10){$L_{k,k-1} $}}
\put(75,50){\makebox(10,10){$L_{2,1} $}}
\put(105,50){\makebox(10,10){$L_{1,0} $}}

\put(15,25){$L_{k+1,k+1-r_2}$}

\put(13,32){\line(-1,1){5}} \put(0,65){\makebox(10,10){$L_{k,k}$}}
\put(30,65){\makebox(10,10){$L_{k-1,k-1} $}}
\put(90,65){\makebox(10,10){$L_{1,1} $}}
\put(120,65){\makebox(10,10){$L_{0,0} $}}
 \put(0,35){\makebox(10,10){$...$}}
\put(90,35){\makebox(10,10){$...$}}
\put(45,20){\makebox(10,10){$...$}}

\put(8,67){\line(1,-1){5}}
\put(32,67){\line(-1,-1){5}} \put(38,67){\line(1,-1){5}}
\put(92,67){\line(-1,-1){5}} \put(98,67){\line(1,-1){5}}
\put(122,67){\line(-1,-1){5}} \put(2,43){\line(-1,1){5}}
\put(8,43){\line(1,1){5}} \put(32,43){\line(-1,1){5}}
\put(38,43){\line(1,1){5}} \put(68,43){\line(1,1){5}}
\put(92,43){\line(-1,1){5}} \put(98,43){\line(1,1){5}}
\put(32,37){\line(-1,-1){9}} \put(38,37){\line(1,-1){9}}
\put(58,33){\line(-1,-1){5}} \put(38,37){\line(1,-1){9}}
\put(72,33){\line(1,-1){5}} \put(92,37){\line(-1,-1){5}}
\put(76,25){$L_{r_2,0}$}

\put(55,70){.}
\put(59,70){.} \put(63,70){.} \put(67,70){.}

 \end{picture}

Case 2a: Diagram for ${I^\alpha(\sigma)}$ when $n+m$ is even and $p$ is odd; $\sigma \leq -{\textstyle \frac{1}{2}}$; $r_2=-\sigma +{\textstyle \frac{1}{2}} \leq [\frac{n+1}{2}]$

\end{center}}

{\tiny
 \vspace{0.2in}
\begin{center}
\setlength{\unitlength}{0.8mm}

\begin{picture}(160,60)(-20,40)
%\begin{picture}(160,140)(-20,-10)

\put(30,95){\makebox(10,10){$L_{k-r_1,k} $}}
\put(60,95){\makebox(10,10){$...$}}
 \put(-15,50){\makebox(10,10){$L_{k+1,k} $}}
\put(-7,60){\line(1,1){5}}
 \put(15,50){\makebox(10,10){$L_{k,k-1} $}}
\put(75,50){\makebox(10,10){$L_{2,1} $}}
\put(105,50){\makebox(10,10){$L_{1,0} $}}

\put(0,65){\makebox(10,10){$L_{k,k}$}}
\put(30,65){\makebox(10,10){$L_{k-1,k-1} $}}
\put(90,65){\makebox(10,10){$L_{1,1} $}}

\put(113,82){.} \put(110,85){.} \put(107,88){.}

\put(120,65){\makebox(10,10){$L_{0,0} $}}

\put(91,99){$L_{0,r_1}$}
\put(32,97){\line(-1,-1){5}}
\put(38,97){\line(1,-1){5}} \put(62,97){\line(-1,-1){5}}
\put(68,97){\line(1,-1){5}} \put(92,97){\line(-1,-1){5}}
\put(98,97){\line(1,-1){5}} \put(8,75){\line(1,1){5}}
\put(32,75){\line(-1,1){5}} \put(38,75){\line(1,1){5}}
\put(92,75){\line(-1,1){5}} \put(98,75){\line(1,1){5}}
\put(122,75){\line(-1,1){5}} \put(8,67){\line(1,-1){5}}
\put(32,67){\line(-1,-1){5}} \put(38,67){\line(1,-1){5}}
\put(92,67){\line(-1,-1){5}} \put(98,67){\line(1,-1){5}}
\put(122,67){\line(-1,-1){5}}

\put(17,82){.} \put(20,85){.} \put(23,88){.}
\put(55,70){.} \put(59,70){.} \put(63,70){.} \put(67,70){.}

 \end{picture}

Case 2a: Diagram for ${I^\alpha(\sigma)}$ when $n+m$ is even and $p$ is odd; $\sigma \geq {\textstyle \frac{1}{2}}$; $r_1=\sigma -{\textstyle \frac{1}{2}} \leq [\frac{n}{2}]$

\end{center}}

\vspace{0.2in}\begin{thm} \label{Ieo}
{\rm (Case 2a)}
Assume that $n+m$ is even, and $p$ is odd.
\begin{enumerate}
\item[(a)] If $\sigma\leq -{\textstyle \frac{1}{2}}$, then
\[\Omega ^{p,q}=L_{\frac{n+1-q}{2},\frac{p-1}{2}}.\]
\vspace{0.2in}
\item[(b)] If $\sigma\geq {\textstyle \frac{1}{2}}$,
then
\[\Omega ^{p,q}=\prec L_{s,t}\succ,\]
where $s =\max\left(0, \frac{n+1-q}{2}\right)$, and $t =\min\left([\frac{n}{2}], \frac{p-1}{2}\right)$.
\end{enumerate}
\end{thm}

\vspace{0.2in}\noindent{\bf Remarks:} \begin{enumerate} \item[(a)] For $(i,j)=(\frac{n+1-q}{2},\frac{p-1}{2})$, we have $i-j=\frac{n+2-m}{2}=-\sigma +{\textstyle \frac{1}{2}}=r_2$. Thus $L_{i,j}$ is at the bottom layer of the module diagram of ${I^\alpha(\sigma)}$, namely is an irreducible submodule. The collection of $L_{i,j}$'s exhausts the set of all irreducible submodules of ${I^\alpha(\sigma)}$.
\item[(b)] If $q\leq n+1$, and $\frac{p-1}{2}\leq [\frac{n}{2}]$, then $s =\max\left(0, \frac{n+1-q}{2}\right)=\frac{n+1-q}{2}$, $t =\min\left([\frac{n}{2}], \frac{p-1}{2}\right)=\frac{p-1}{2}$, and $t-s=\frac{m-(n+2)}{2}=\sigma -{\textstyle \frac{1}{2}} =r_1$. Such $L_{s,t}$'s are exactly those at the top layer of the module diagram of ${I^\alpha(\sigma)}$, namely irreducible quotients. The rest of $L_{s,t}$'s are those on the ``left boundary" ($t=k=[\frac{n}{2}]$) or the ``right boundary" ($s=0$). Note that when $n$ is even, the submodule $L_{k+1,k}$ is empty.
\end{enumerate}

\vspace{0.2in}

Next recall the module diagram of ${I^\alpha(\sigma)}$ for Case 2b:

{\tiny
 \vspace{0.2in}
\begin{center}
\setlength{\unitlength}{0.8mm}

\begin{picture}(160,60)(-20,10)
%\begin{picture}(160,140)(-20,-10)

\put(-15,50){\makebox(10,10){$L_{0,0} $}}
\put(-7,60){\line(1,1){5}}
\put(15,50){\makebox(10,10){$L_{1,1} $}}
\put(75,50){\makebox(10,10){$L_{k-1,k-1} $}}
\put(105,50){\makebox(10,10){$L_{k,k} $}}
\put(15,25){$L_{0,r_1}$}
\put(13,32){\line(-1,1){5}} \put(0,65){\makebox(10,10){$L_{1,0}$}}
\put(30,65){\makebox(10,10){$L_{2,1} $}}
\put(90,65){\makebox(10,10){$L_{k,k-1} $}}
\put(120,65){\makebox(10,10){$L_{k+1,k} $}}
 \put(0,35){\makebox(10,10){$...$}}
\put(90,35){\makebox(10,10){$...$}}
\put(45,20){\makebox(10,10){$...$}}
\put(8,67){\line(1,-1){5}}
\put(32,67){\line(-1,-1){5}} \put(38,67){\line(1,-1){5}}
\put(92,67){\line(-1,-1){5}} \put(98,67){\line(1,-1){5}}
\put(122,67){\line(-1,-1){5}} \put(2,43){\line(-1,1){5}}
\put(8,43){\line(1,1){5}} \put(32,43){\line(-1,1){5}}
\put(38,43){\line(1,1){5}} \put(68,43){\line(1,1){5}}
\put(92,43){\line(-1,1){5}} \put(98,43){\line(1,1){5}}
\put(32,37){\line(-1,-1){9}} \put(38,37){\line(1,-1){9}}
\put(58,33){\line(-1,-1){5}} \put(38,37){\line(1,-1){9}}
\put(72,33){\line(1,-1){5}} \put(92,37){\line(-1,-1){5}}
\put(76,25){$L_{k-r_1,k}$}
\put(55,70){.}
\put(59,70){.} \put(63,70){.} \put(67,70){.}

 \end{picture}

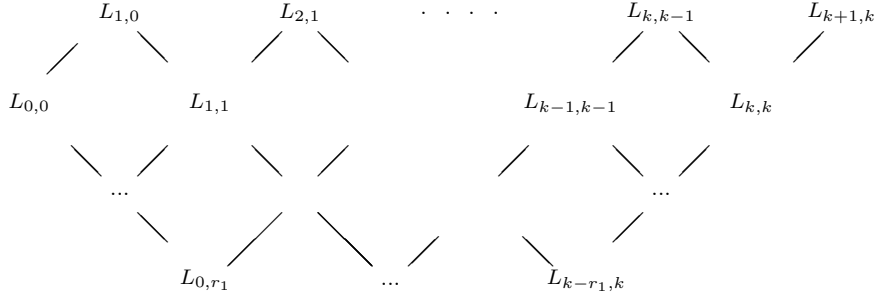
Case 2b: Diagram for ${I^\alpha(\sigma)}$ when $n+m$ and $p$ are even; $\sigma \leq -{\textstyle \frac{1}{2}}$; $r_1=-\sigma -{\textstyle \frac{1}{2}} \leq [\frac{n}{2}]$

\end{center}}

{\tiny
 \vspace{0.2in}
\begin{center}
\setlength{\unitlength}{0.8mm}

\begin{picture}(160,60)(-20,40)
%\begin{picture}(160,140)(-20,-10)

\put(30,95){\makebox(10,10){$L_{r_2,0} $}}
\put(60,95){\makebox(10,10){$...$}}
 \put(-15,50){\makebox(10,10){$L_{0,0} $}}
\put(-7,60){\line(1,1){5}}
 \put(15,50){\makebox(10,10){$L_{1,1} $}}
\put(75,50){\makebox(10,10){$L_{k-1,k-1} $}}
\put(105,50){\makebox(10,10){$L_{k,k} $}}
\put(0,65){\makebox(10,10){$L_{1,0}$}}
\put(30,65){\makebox(10,10){$L_{2,1} $}}
\put(90,65){\makebox(10,10){$L_{k,k-1} $}}

\put(113,82){.} \put(110,85){.} \put(107,88){.}

\put(120,65){\makebox(10,10){$L_{k+1,k} $}}
\put(91,99){$L_{k+1,k+1-r_2}$}
\put(32,97){\line(-1,-1){5}}
\put(38,97){\line(1,-1){5}} \put(62,97){\line(-1,-1){5}}
\put(68,97){\line(1,-1){5}} \put(92,97){\line(-1,-1){5}}
\put(98,97){\line(1,-1){5}} \put(8,75){\line(1,1){5}}
\put(32,75){\line(-1,1){5}} \put(38,75){\line(1,1){5}}
\put(92,75){\line(-1,1){5}} \put(98,75){\line(1,1){5}}
\put(122,75){\line(-1,1){5}} \put(8,67){\line(1,-1){5}}
\put(32,67){\line(-1,-1){5}} \put(38,67){\line(1,-1){5}}
\put(92,67){\line(-1,-1){5}} \put(98,67){\line(1,-1){5}}
\put(122,67){\line(-1,-1){5}}
\put(17,82){.} \put(20,85){.} \put(23,88){.}
 \put(55,70){.}
\put(59,70){.} \put(63,70){.} \put(67,70){.}

 \end{picture}

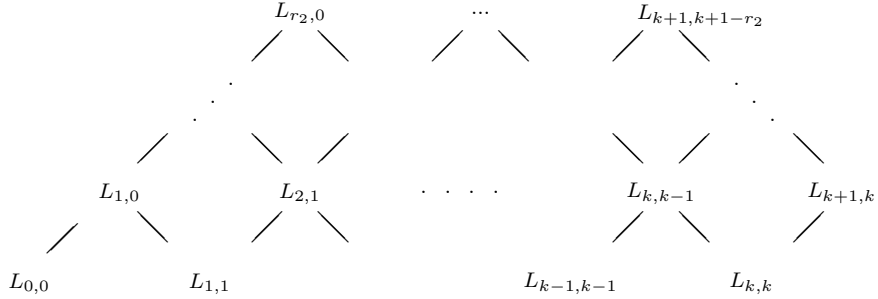
Case 2b: Diagram for ${I^\alpha(\sigma)}$ when both $n+m$ and $p$ are even; $\sigma \geq {\textstyle \frac{1}{2}}$; $r_2=\sigma +{\textstyle \frac{1}{2}} \leq [\frac{n+1}{2}]$

\end{center}}

\vspace{0.2in}\begin{thm} \label{Iee}
{\rm (Case 2b)}
Assume that $n+m$ is even, and $p$ is even.
\begin{enumerate}
\item[(a)] If $\sigma\leq -{\textstyle \frac{1}{2}}$, then
\[\Omega ^{p,q}=L_{\frac{p}{2},\frac{n-q}{2}}.\]
\vspace{0.2in}
\item[(b)] If $\sigma\geq {\textstyle \frac{1}{2}}$,
then
\[\Omega ^{p,q}=\prec L_{s,t}\succ,\]
where $s =\min\left([\frac{n+1}{2}], \frac{p}{2}\right)$, and $t =\max\left(0, \frac{n-q}{2}\right)$.
\end{enumerate}
\end{thm}

\vspace{0.2in}\noindent{\bf Remarks:} \begin{enumerate} \item[(a)] For $(i,j)=(\frac{p}{2},\frac{n-q}{2})$, we have $j-i=\frac{n-m}{2}=-\sigma -{\textstyle \frac{1}{2}}=r_1$. Thus $L_{i,j}$ is at the bottom layer of the module diagram of ${I^\alpha(\sigma)}$, namely is an irreducible submodule. The collection of $L_{i,j}$'s exhausts the set of all irreducible submodules of ${I^\alpha(\sigma)}$.
\item[(b)] If $\frac{p}{2}\leq [\frac{n+1}{2}]$ and $q\leq n$, then $s =\min\left([\frac{n+1}{2}], \frac{p}{2}\right)=\frac{p}{2}$, $t =\max\left(0, \frac{n-q}{2}\right)= \frac{n-q}{2}$, and $s-t=\frac{m-n}{2}=\sigma +{\textstyle \frac{1}{2}} =r_2$. Such $L_{s,t}$'s are exactly those at the top layer of the module diagram of ${I^\alpha(\sigma)}$, namely irreducible quotients. The rest of $L_{s,t}$'s are those on the ``left boundary" ($t=0$) or the ``right boundary" ($s=[\frac{n+1}{2}]$). Note that when $n$ is even, the subquotient $L_{k+1,j}$ is empty.
\end{enumerate}

\end{document}